\documentclass[12pt,fleqn]{article}
\usepackage{greg}
\usepackage{titlesec}
\usepackage{titling}
\usepackage[nottoc]{tocbibind}
\usepackage{fancyvrb}
\usepackage[font=footnotesize,labelfont=bf]{caption}
\numberwithin{equation}{section}
\graphicspath{ {./img/} }

\renewenvironment{abstract}{\section*{\abstractname}}{}
\titleformat*{\section}{\large\bfseries\itshape}
\titleformat*{\subsection}{\normalsize\bfseries}
\titleformat*{\subsubsection}{\normalsize\sffamily}
\titleformat{\chapter}[hang]{\LARGE\sffamily}{\LARGE\thechapter}{1ex}{}[]
\titleformat{name=\chapter,numberless}[hang]{\LARGE\bfseries}{}{0ex}{}[]

\title{Discrete differential geometry in homotopy type theory}
\author{Greg Langmead}
\begin{document}
\begin{titlepage}
\maketitle
\thispagestyle{empty}

\begin{center}
Department of Philosophy\\
Carnegie Mellon University\\
Pittsburgh, PA 15213
\vspace{1in}

Submitted in partial fulfillment of the requirements for the degree of\\
Master of Science in Logic, Computation and Methodology
\vspace{1in}

\textbf{Thesis Committee}

Steven M. Awodey, chair\\
Mathieu Anel

\end{center}
\end{titlepage}

\thispagestyle{empty}
\begin{dedication}
This thesis is dedicated to John Baez, Sean M. Carroll, Sabine Hossenfelder, and other communicators who are carrying the torch of science forward in the spirit of my hero Carl Sagan. I have followed you all for many years, and you have inspired me to continue my studies alongside my career. Thank you.
\end{dedication}

\clearpage

\thispagestyle{empty}
\begin{abstract}
Type families on higher inductive types such as pushouts can capture homotopical properties of differential geometric constructions including connections, curvature, and vector fields. We define a class of pushouts based on simplicial complexes, then define principal bundles, connections, and curvature on these. We provide an example of a tangent bundle but do not prove when these must exist. We define vector fields, and the index of a vector field. Our main result is a theorem relating total curvature and total index, a key step to proving the Gauss-Bonnet theorem and the Poincaré-Hopf theorem, but without an existing definition of Euler characteristic to compare them to. We draw inspiration in part from the young field of discrete differential geometry, and in part from the original classical proofs, which often make use of triangulations and other discrete arguments.
\end{abstract}
\clearpage

\newlength{\mylen}
\newlength{\mylin}

\thispagestyle{empty}
\tableofcontents 
\clearpage
\setcounter{page}{5}
\pagestyle{plain}

\begin{quote} 
``It is always ourselves we work on, whether we realize it or not. There is no other work to be done in the world.''
 --- Stephen Talbott, \emph{The Future Does Not Compute}\cite{talbott}
\end{quote}

\section{Overview}

We will define 
\begin{itemize}
\item principal bundles in Section~\ref{sec:torsors},
\item simplicial complexes, and homotopical realizations of these in Section~\ref{sec:discrete_man},
\item tangent bundles and vector fields in Section~\ref{sec:vector_fields},
\end{itemize}
and observe emerging from those definitions the presence of
\begin{itemize}
\item connections and curvature in Section~\ref{sec:connections},
\item the index of a vector field in Section~\ref{sec:totals},
\end{itemize}
and then define in Section~\ref{sec:totals}, for a 2-dimensional simplicial complex
\begin{itemize}
\item the total curvature, as in the Gauss-Bonnet theorem
\item the total index of a vector field, as in the Poincaré-Hopf theorem,
\item and prove the equality of these to each other when the complex is oriented (Theorem~\ref{thm:total_index_total_curvature}).
\end{itemize}

We will build up an example of all of these structures on an octahedron model of the sphere. We will not, however, be supplying a separate definition of Euler characteristic so as to truly reproduce the Gauss-Bonnet and Poincaré-Hopf theorems.

Once we have defined homotopical realizations of simplicial complexes in Section~\ref{sec:discrete_man}, we will focus on dimensions 1 and 2. In dimension 1 we define polygons, which we prove are equivalent to \( S^1 \), and so give terms in the type \( \EMzo\defeq \sit{A:\uni}||A\simeq S^1||_1 \). We can call this component of the universe ``mere circles.'' In dimension 2 we will focus on a subset of complexes where the neighboring vertices and edges of each vertex (the vertex's ``link'') form a polygon. The homotopical realization \( \mm \) of such a complex then has a map \( \link \) from each vertex to a homotopical polygon, i.e. a map to \( \EMzo \). 

Given a map \( \mm\to\EMzo \) we can form the pullback
\begin{center}
\begin{tikzcd}[cramped]
  {P} & {\mathrm{EM}_\bullet(\mathbb{Z},1)} \\
  \mm & {\mathrm{EM}(\mathbb{Z},1)}
  \arrow["", from=1-1, to=1-2]
  \arrow["{\mathrm{pr}_1}"', from=1-1, to=2-1]
  \arrow["\lrcorner"{anchor=center, pos=0.125}, draw=none, from=1-1, to=2-2]
  \arrow["{\mathrm{pr}_1}"', from=1-2, to=2-2]
  \arrow["{\mathsf{link}}", from=2-1, to=2-2]
\end{tikzcd}
\end{center}
to obtain a bundle of mere circles. We will discuss when such a map factors through a map \( \mm\to \K(\zz,2) \), and hence when it is in fact a principal fibration.

Then in Section~\ref{sec:connections} we will name various elements of the above construction, indicating their relationship to classical definitions.

In Section~\ref{sec:vector_fields} we will define vector fields, which require a tangent bundle. We will introduce a method for decomposing a vector field along a concatenation of paths.

Finally, in Section~\ref{sec:totals} we will define a method for visiting all the faces of a manifold in order to form ``totals'' of local objects. We will examine the total curvature and the total index and prove that they are equal. Our proof tracks very closely with the classical proof of Hopf\cite{hopf}, presented in detail in Needham\cite{needham}.

This work is appearing while other researchers are actively exploring and exploiting the related notion of CW complexes in HoTT, for example \cite{ljungstrom_cellular}. Combining and organizing these approaches would likely provide additional visibility into a program for importing more of 20th century topology and gauge theory into HoTT.

\clearpage
\section{Torsors and principal bundles}
\label{sec:torsors}

\subsection{Torsors}
\label{subsec:torsors}
We will review some definitions and facts, drawing on the excellent resource \cite{buchholtz2023central}. Let \( G \) be a group, possibly a higher group.
\begin{mydef}
A \defemph{right \( G \)-object} is a type \( X \) equipped with a homomorphism \( \phi:G^{\mathrm{op}}\to\Aut(X) \). If in addition the map \( (\pr_1,\phi):X\times G\to X\times X \) is an equivalence then we say \( (X,\phi) \) is a \defemph{\( G \)-torsor}. Denote the type of \( G \)-torsors by \( BG \). Denote the \( G \)-torsor given by \( G \) itself under right-multiplication by \( \reg{G} \).
\end{mydef}

A \( G \)-equivariant map from \( (X,\phi) \) to \( (Y,\psi) \) is a function \( f:X\to Y \) such that \( f(\phi(g,x))=\psi(g,f(x)) \). Denote the type of \( G \)-equivariant maps by \( X\to_G Y \).

\begin{mylemma}
(\cite{buchholtz2023central} Lemma 5.2). If \( (X,\phi),(Y,\psi):BG \) then there is a natural equivalence \( (X=_{BG}Y) \simeq (X\to_G Y) \).\qed
\end{mylemma}

\begin{mylemma}
(\cite{buchholtz2023central} Proposition 5.4). A \( G \)-object \( (X,\phi) \) is a \( G \)-torsor if and only if there merely exists a \( G \)-equivariant equivalence \( \reg{G}\to_G X \).\qed
\end{mylemma}

These lemmas lead to

\begin{mycor}
(\cite{buchholtz2023central} Corollary 5.5). The pointed type \( (BG,\reg{G}) \) is a \( \K(G,1) \), i.e. \( BG \) is connected and has a pointed equivalence \( \loopy BG\simeq_* G \).\qed
\end{mycor}

We call types \( X \) such that \( \loopy X\simeq_* G \) a \emph{delooping} of \( G \). So the type of \( G \)-torsors \( BG \) is a delooping of \( G \).

If \( G \) is abelian then \( \K(G,1) \) can be delooped, and in fact can be delooped in a unique way, any number of times. In the next section we will learn how to proceed from a construction of \( \K(G, 1) \) to obtain \( \K(G, n) \).

\subsection{Bundles of Eilenberg-Mac Lane spaces}

To construct maps into \( \Kzt \) we will follow Scoccola\cite{sco}. When can a map into a connected component of the universe be factored through an Eilenberg-Mac Lane space?

\begin{mydef}
Let \( G \) be a group. Let \[ \EM(G,n)\defeq \BAut(\K(G,n))\defeq \sit{Y:\uni}||Y\simeq \K(G,n)||_{-1}\] be the connected component of \( \uni \) containing \( \K(G,n) \). A \defemph{\( \K(G,n) \)-bundle} on a type \( M \) is a map \( M\to\EM(G,n) \).
\end{mydef}

Scoccola constructs a map by composing 
\begin{itemize}
\item suspension \( \susp:\EM(G,n)\to\EM_{\bullet\bullet}(G,n) \) (see \cite{hottbook} §6.5), which maps into types with two points (denoted by the bullets),
\item \( (n+1) \)-truncation (see \cite{hottbook} §7.3),
\item forgetting a point \( F_\bullet:\EM_{\bullet\bullet}(G,n)\to \EM_{\bullet}(G,n) \),
\end{itemize}
to form the composition
\[ 
\EM(G,n)\xrightarrow[]{||\susp||_{n+1}} \EM_{\bullet\bullet}(G,n+1)\xrightarrow[]{F_\bullet}\EMp(G,n+1)
\]
to types with two points (north and south), then to pointed types (by forgetting the south point).

\begin{mydef}
Given \( f:M\to\EM(G,n) \), the \defemph{associated action of \( M \) on \( G \)}, denoted by \( f_\bullet \) is defined to be \( f_\bullet=F_\bullet\circ||\susp||_{n+1}\circ f \).
\end{mydef}

\begin{mythm}
\label{thm:sco}
(Scoccola\cite{sco} Proposition 2.39). A \( \K(G,n) \) bundle \( f:M\to\EM(G,n) \) factors through a map \( M\to\K(G,n+1) \), and so is a principal fibration, if and only if the associated action \( f_\bullet \) is merely homotopic to a constant map.
\end{mythm}

The above theory of classifying spaces should be relevant to any future project to bring the study of gauge theory into homotopy type theory. In this note we will be focused on the special case \( G=\zz \) and \( \EMzo \) and \( \Kzt \).

\begin{mynote}
Iterating the loop map gives the isomorphism \( \loopy^{(n+1)}:\EMp(G,n+1)\simeq \K(\Aut G,1) \) (see \cite{sco} Lemma 2.7). Theorem~\ref{thm:sco} therefore says that the map \( f \) factors through \( \K(G,n+1) \) if and only if the map into \( \K(\Aut G,1) \) is homotopic to a constant. In the case of \( G=\zz \), the map \( f_\bullet:M\to \K(\Aut \zz, 1) \) deserves to be called the first Stiefel-Whitney class of \( f \) and one can interpret its triviality as \emph{orientability}. . This point of view is discussed in Schreiber\cite{dcct} (starting with Example 1.2.138) and in Myers\cite{myersgood}.
\end{mynote}

\subsection{Pathovers}
\label{sec:pathovers}
Here we will recall basic facts from homotopy type theory. See for example \cite{hottbook} or \cite{egbert}. Suppose we have \( T:M\to\uni \) and \( P\defeq\sit{x:M}Tx \). Recall that if \( p:a=_M b \) then \( T \) acts on \( p \) with what's called the \emph{action on paths}, denoted \( \ap(T)(p):Ta=Tb \) (\cite{hottbook} §2.2). This is a path in the codomain \( \uni \) of \( T \). Type theory also provides a function called \emph{transport}, denoted \( \tr(p):Ta\to Tb \) (\cite{hottbook} §2.3) which acts on the fibers of \( P \). \( \tr(p) \) is a function, acting on the terms of the types \( Ta \) and \( Tb \), and univalence tells us this is the isomorphism corresponding to \( \ap(T)(p) \).

\begin{figure}[H]
\centering
\includegraphics[width=200pt]{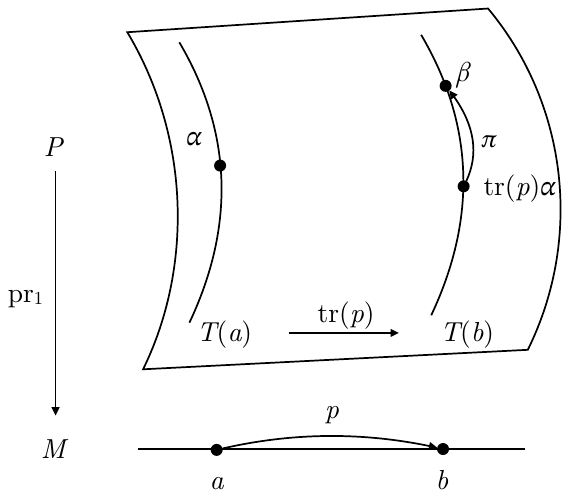}
\caption{A path \( \pi \) over the path \( p \) in the base involves the transport function.}
\label{fig:pathovers}
\end{figure}

Type theory also tells us that paths in \( P \) are given by pairs of paths: a path \( p:a=_M b \) in the base, and a pathover \( \pi:\tr(p)(\alpha)=_{Tb}\beta \) between \( \alpha:Ta \) and \( \beta:Tb \) in the fibers (\cite{hottbook} §2.7). We can't directly compare \( \alpha \) and \( \beta \) since they are of different types, so we apply transport to one of them. We say \( \pi \) lies over \( p \). See Figure~\ref{fig:pathovers}.

\begin{figure}[h]
\centering
\begin{tikzpicture}[
node distance = 20mm and 20mm,
V/.style = {circle, fill, draw=black, inner sep=1pt},
every edge quotes/.style = {auto},
arrow/.style={->,semithick}
]
\begin{scope}[nodes=V]
  \node[label=above left:\( \phi(a) \)] (1) {};
  \node[label=above right:\( \phi(a') \)] (2) [right=of 1]  {};
  \node[label=below right:\( \psi(a') \)] (3) [below=of 2]  {};
  \node[label=below left:\( \psi(a) \)] (4) [below=of 1]  {};
  \node[label=below:\( a \)] (5) [below=of 4]  {};
  \node[label=below:\( a' \)] (6) [below=of 3]  {};
\end{scope}
\draw[arrow]
        (2)  edge[swap, "\( \phi(p)^{-1} \)"] (1)
        (4)  edge["\( \psi(p) \)"] (3)
        (1)  edge[swap, "\( q \)"] (4)
        (5)  edge["\( p \)"] (6);
\end{tikzpicture}
\caption{Transport along \( p \) in the fibers of a family of paths. The fiber over \( a \) is \( \phi(a)=\psi(a) \) where \( \phi,\psi:A\to B \).}
\label{fig:transport_family_of_paths}
\end{figure}
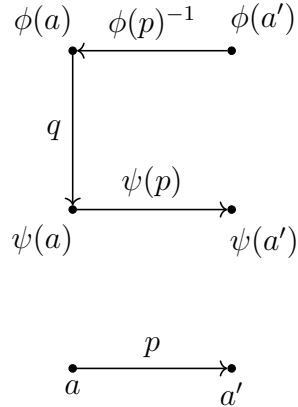

Given functions \( \phi,\psi:A\to B \) between two arbitrary types we can form a type family of paths \( \alpha:A\to\uni \) by \( \alpha(a)\defeq(\phi(a)=_B\psi(a)) \). Transport in this family is given by concatenation as follows (see Figure~\ref{fig:transport_family_of_paths}), where \( p:a=_A a' \) and \( q:\phi(a)=\psi(a) \) (\cite{hottbook} Theorem 2.11.3):
\[ 
\tr(p)(q) = \phi(p)^{-1}\cdot q\cdot \psi(p)
\]
which gives a path in \( \phi(a')=\psi(a') \) by connecting dots between the terms \( \phi(a'), \phi(a), \psi(a), \psi(a') \). This relates a would-be homotopy \( \phi\sim\psi \) specified at a single point, to a point at the end of a path. We will use this to help construct such homotopies.

Finally, recall that in the presence of a section \( X:M\to P \) there is a dependent generalization of \( \ap \) called \( \apd \): \( \apd(X)(p):\tr(p)(X(a))=X(b) \) which is a pathover between the two values of the section over the basepoints of the path \( p \) (\cite{hottbook} Lemma 2.3.4).

\clearpage
\section{Discrete manifolds}
\label{sec:discrete_man}
We will remind ourselves of the definition of a classical simplicial complex, in sets. Then we will create a type that realizes the data of a complex, using pushouts.

\subsection{Abstract simplicial complexes}

\begin{mydef}
An \defemph{abstract simplicial complex \( M \) of dimension \( n \)} is an ordered list \( M\defeq[M_0,\ldots,M_n] \) consisting of a set \( M_0 \) of vertices, and for each \( 0<k\leq n \) a set \( M_k \) of subsets of \( M_0 \) of cardinality \( k+1 \), such that for any \( j<k \), any \( (j+1) \)-element subset of an element of \( M_k \) is an element of \( M_j \). The elements of \( M_k \) are called \defemph{\( k \)-faces}. A \defemph{morphism} \( f \) from \( M=[M_0,\ldots,M_m] \) to \( N=[N_0,\ldots,N_n] \) is a function on vertices \( f:M_0\to N_0 \) such that for any face of \( M \) the image under \( f \) is a face of \( N \). Denote by \( \simcomp \) the type of abstract simplicial complexes of dimension \( n \). Let \( M_{\leq k}\defeq [M_0,\ldots,M_k] \). We call \( M_{\leq k} \) the \defemph{\( k \)-skeleton} of \( M \), and it is a (\( k \)-)complex in its own right. \( M \) is automatically equipped with a chain of inclusions of the skeleta \( M_0\hookrightarrow M_{\leq 1}\hookrightarrow\cdots\hookrightarrow M_{\leq n}=M \), which are simply the inclusions of lists into longer lists.
\end{mydef}

We can imagine constructing a simplicial complex dimension by dimension with a simple procedure of taking unions. This will parallel the realization construction that follows.

\begin{mydef}
Denote by \( \Delta(n) \) the simplicial complex obtained by taking the set of all subsets of the standard \( (n+1) \)-element set \( \Delta(n)_0\defeq\{0,\ldots,n\} \). We call \( \Delta(n) \) the \defemph{complete \( n \)-simplex}. We will refer to the \( (n-1) \)-skeleton \( \Delta(n)_{\leq(n-1)} \) with the suggestive notation \( \partial \Delta(n) \).
\end{mydef}

Note that \( \Delta(n)_{(n-1)} \) has \( n+1 \) elements. For example, \( \Delta(2)_1 \) consists of the three 2-element subsets of \( \{0, 1, 2\} \). 

Given \( M=[M_0,\ldots,M_n] \) and \( k\leq n \), a face \( f_k:M_k \) is the union of \( (k+1) \) faces in \( M_{k-1} \), and so the \( k \)-skeleton is obtained from the \( (k-1) \)-skeleton by forming the following pushout of sets:
\begin{equation}
\begin{tikzcd}
  {M_k\times\partial \Delta(k)} & {M_k} \\
  {M_{\leq (k-1)}} & {M_{\leq k}}
  \arrow["{\mathrm{pr}_1}", from=1-1, to=1-2]
  \arrow["{a_k}"', from=1-1, to=2-1]
  \arrow[from=1-2, to=2-2]
  \arrow[from=2-1, to=2-2]
  \arrow["\ulcorner"{anchor=center, pos=0.125, rotate=180}, draw=none, from=2-2, to=1-1]
\label{eq:attach}
\end{tikzcd}
\end{equation}

where the vertical ``attach'' map \( a_k(f_k, -) \) picks out the \( k+1 \) subsets of \( f_k \).

\begin{mydef}
In an abstract simplicial complex \( M \) of dimension \( n \), the \defemph{link} of a vertex \( v \) is the \( (n-1) \)-dimensional subcomplex containing every face \( m\in M_{n-1} \) such that \( v\notin m \) and \( m\cup v \) is an \( n \)-face of \( M \).
\label{def:link}
\end{mydef}

\begin{figure}[h]
\centering
\begin{tikzpicture}
  \draw
    (0, 0) grid[step=1cm] (3, 3)
    (0, 2) edge[ultra thick] (1, 3)
    (0, 1) edge[ultra thick] (0, 2)
    (0, 1) edge[ultra thick] (1, 1)
    (1, 1) edge[ultra thick] (2, 2)
    (2, 2) edge[ultra thick] (2, 3)
    (1, 3) edge[ultra thick] (2, 3)
    (0, 1) -- (2, 3)
    (0, 0) -- (3, 3)
    (1, 0) -- (3, 2)
    (2, 0) -- (3, 1);
  \fill[radius=3pt]
    \foreach \x in {0, ..., 3} {
      \foreach \y in {0, ..., 3} {
        (\x, \y) circle[]
    }};
  \path[above left]
    \foreach \p/\v in {
      {1, 3}/a,
      {2, 3}/b,
      {0, 2}/f,
      {1, 2}/v,
      {2, 2}/c,
      {0, 1}/e,
      {1, 1}/d%
    } {
      (\p) node {$\v$}
    };
\end{tikzpicture}
\caption{The link of \( v \) in this complex consists of the vertices \( \{a,b,c,d,e,f\} \) and the edges \( \{ab,bc,cd,de,ef,fa\} \), forming a hexagon.}
\label{fig:link}
\end{figure}
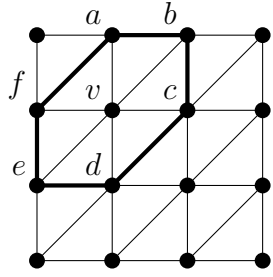

The link is easier to understand as all the neighboring vertices of \( v \) and the subcomplex containing these. See for example Figure~\ref{fig:link}.

\begin{mynote}
The \emph{geometric realization} of a complex uses the combinatorial data to form pushouts of standard simplices inside the category of topological spaces. A \emph{simplicial sphere} is a simplicial complex whose geometric realization is homeomorphic to a sphere. A classical 1940 result of Whitehead, building on Cairn, states that every smooth \( n \)-manifold is the geometric realization of a simplicial complex of dimension \( n \) such that the link is the geometric realization of an \( (n-1) \)-sphere\cite{whitehead_triangulation}. For more of this theory see the classic book by Kirby and Siebenmann\cite{kirby_siebenmann}.
\end{mynote}

\subsection{Higher inductive realizations}

We will realize a simplicial complex as a higher inductive type by forming a sequence of pushouts. We will work upward by dimension so as to define a standard sphere that we need in the next dimension.

\begin{mydef}
The \defemph{realization \( \mm \) of a 0-dimensional simplicial complex} \( M \) is simply the set \( M_0 \).
\end{mydef}

\begin{mydef}
The \defemph{simplicial 0-sphere} \( \bdsimplexn{1} \) is the realization of \( \partial \Delta(1) \).
\end{mydef}

\begin{mydef}
The \defemph{realization \( \mm \) of a 1-dimensional simplicial complex} \( M=[M_0, M_1] \) is the pushout
\end{mydef}
\begin{center}
\begin{tikzcd}
  {M_1\times \bdsimplexn{1}} & {M_1} \\
  {M_0=\mathbb{M}_0} & {\mathbb{M}_1}
  \arrow["{\mathrm{pr}_1}", from=1-1, to=1-2]
  \arrow["{\mathbb{A}_0}"', from=1-1, to=2-1]
  \arrow["{*_{\mathbb{M}_1}}", from=1-2, to=2-2]
  \arrow["{h_1}"', shorten <=14pt, shorten >=14pt, Rightarrow, from=2-1, to=1-2]
  \arrow[from=2-1, to=2-2]
  \arrow["\ulcorner"{pos=0, rotate=180}, shift left=1, draw=none, from=2-2, to=1-1]
\end{tikzcd}
\end{center}
where \( a_k \) is the attachment data of the edges in \( M_1 \), as in the diagram~\ref{eq:attach}. The right vertical map \( *_{\mm_1} \) provides a hub point for each edge, and the homotopy \( h_1 \) provides the spokes that connect the hub to the vertices.

\begin{mydef}
The \defemph{simplicial 1-sphere} \( \bdsimplexn{2} \) is the realization of \( \partial \Delta(2) \).
\end{mydef}
\begin{center}
\begin{tikzcd}
  {\partial \Delta(2)\times \bdsimplexn{1}} & {\partial \Delta(2)} \\
  {\Delta(2)_0} & {\bdsimplex{2}}
  \arrow["{\mathrm{pr}_1}", from=1-1, to=1-2]
  \arrow["{\mathbb{A}_0}"', from=1-1, to=2-1]
  \arrow["{*_{\bdsimplexn{2}}}", from=1-2, to=2-2]
  \arrow["{h_1}"', shorten <=11pt, shorten >=11pt, Rightarrow, from=2-1, to=1-2]
  \arrow[from=2-1, to=2-2]
  \arrow["\ulcorner"{anchor=center, pos=0.125, rotate=180}, draw=none, from=2-2, to=1-1]
\end{tikzcd}
\end{center}

\begin{mydef}
A \defemph{realization \( \mm \) of a 2-dimensional simplicial complex} \( [M_0, M_1, M_2] \) is the pushout
\end{mydef}
\begin{center}
\begin{tikzcd}
  & {M_2\times \bdsimplexn{2}} & {M_2} \\
  {\mathbb{M}_0} & {\mathbb{M}_1} & {\mathbb{M}_2} \\
  {M_1\times \bdsimplexn{1}} & {M_1}
  \arrow["{\mathrm{pr}_1}", from=1-2, to=1-3]
  \arrow["{\mathbb{A}_1}"', from=1-2, to=2-2]
  \arrow["{*_{\mathbb{M}_2}}", from=1-3, to=2-3]
  \arrow[from=2-1, to=2-2]
  \arrow["{h_1}"', shorten <=27pt, shorten >=27pt, Rightarrow, from=2-1, to=3-2]
  \arrow["{h_2}", shorten <=17pt, shorten >=17pt, Rightarrow, from=2-2, to=1-3]
  \arrow[from=2-2, to=2-3]
  \arrow["\ulcorner"{pos=-0.1, rotate=-90}, draw=none, from=2-2, to=3-1]
  \arrow["\ulcorner"{pos=0, rotate=180}, draw=none, from=2-3, to=1-2, shift left=1]
  \arrow["{\mathbb{A}_0}", from=3-1, to=2-1]
  \arrow["{\mathrm{pr}_1}"', from=3-1, to=3-2]
  \arrow["{*_{\mathbb{M}_1}}"', from=3-2, to=2-2]
\end{tikzcd}

\end{center}
where \( \aaa_1 \) is such that this additional diagram commutes
\begin{center}
\begin{tikzcd}
  {M_2\times\partial \Delta(2)} & {M_2\times\bdsimplexn{2}} \\
  {M_{\leq 1}} & {\mathbb{M}_1}
  \arrow["{\mathrm{id}\times *_{\bdsimplexn{2}}}", from=1-1, to=1-2]
  \arrow["{a_1}"', from=1-1, to=2-1]
  \arrow["{\mathbb{A}_1}", from=1-2, to=2-2]
  \arrow["{*_{\mathbb{M}_{\leq 1}}}"', from=2-1, to=2-2]
\end{tikzcd}
\end{center}
In this diagram \( a_1 \) is the simplicial attaching map from \ref{eq:attach}, and \( *_{\mathbb{M}_{\leq 1}} \) simply gathers the hub maps from \( M_0 \) and \( M_1 \) into \( \mm_1 \). The commutativity then says that \( \aaa_1 \) reflects the attachment data.
\begin{mydef}
Given a notion of realization in dimension \( n-1 \), a \defemph{realization \( \mm \) of an \( n \)-dimensional simplicial complex} \( M=[M_0, \ldots, M_n] \) is the pushout
\end{mydef}
\begin{center}%
\begin{tikzcd}
  \cdots & {M_{n-2}} & {M_n\times \bdsimplexn{n}} & {M_n} \\
  \cdots & {\mathbb{M}_{n-2}} & {\mathbb{M}_{n-1}} & {\mathbb{M}_n} \\
  & \cdots & {M_{n-1}}
  \arrow[from=1-1, to=1-2]
  \arrow["{*_{\mathbb{M}_{n-2}}}", from=1-2, to=2-2]
  \arrow["{\mathrm{pr}_1}", from=1-3, to=1-4]
  \arrow["{\mathbb{A}_{n-1}}"', from=1-3, to=2-3]
  \arrow["{*_{\mathbb{M}_{n}}}", from=1-4, to=2-4]
  \arrow["{h_{n-2}}", shorten <=8pt, shorten >=8pt, Rightarrow, from=2-1, to=1-2]
  \arrow[from=2-1, to=2-2]
  \arrow["\ulcorner"{anchor=center, pos=0.125, rotate=180}, draw=none, from=2-2, to=1-1]
  \arrow[from=2-2, to=2-3]
  \arrow["{h_{n-1}}"', shorten <=10pt, shorten >=10pt, Rightarrow, from=2-2, to=3-3]
  \arrow["{h_n}", shorten <=10pt, shorten >=10pt, Rightarrow, from=2-3, to=1-4]
  \arrow[from=2-3, to=2-4]
  \arrow["\ulcorner"{anchor=center, pos=0.125, rotate=-90}, draw=none, from=2-3, to=3-2]
  \arrow["\ulcorner"{anchor=center, pos=0.125, rotate=180}, draw=none, from=2-4, to=1-3]
  \arrow[from=3-2, to=3-3]
  \arrow["{*_{\mathbb{M}_{n-1}}}"', from=3-3, to=2-3]
\end{tikzcd}
\end{center}
where \( \aaa_n \) is such that the following commutes
\begin{center}
\begin{tikzcd}
  {M_n\times\partial \Delta(n)} & {M_n\times\bdsimplexn{n}} \\
  {M_{\leq n}} & {\mathbb{M}_n}
  \arrow["{\mathrm{id}\times *_{\bdsimplexn{n}}}", from=1-1, to=1-2]
  \arrow["{a_1}"', from=1-1, to=2-1]
  \arrow["{\mathbb{A}_n}", from=1-2, to=2-2]
  \arrow["{*_{\mathbb{M}_{\leq n}}}"', from=2-1, to=2-2]
\end{tikzcd}
\end{center}
Gathering some of the inductive process into one definition, we have
\begin{mydef}
\label{def:higher_realization} 
A \defemph{realization} \( \mm \) of an abstract simplicial complex \( M:\simcomp \) consists of 
\begin{enumerate}
\item \( n+1 \) types \( \mm_0,\ldots,\mm_n \) where \( \mm_0\defeq M_0 \),
\item \( n \) spans \( \myspan{\mm_{i}}{M_{i+1}\times \bdsimplexn{i+1}}{M_{i+1}}{\aaa_{i}}{\pr_1} \), \( i=0,\ldots,n-1 \), where \( \bdsimplexn{i+1} \) is the realization of the boundary of a standard simplex and \( \aaa_i \) are called \defemph{attachment maps},
\item \( n \) pushout squares from each span to \( \mm_{i+1} \), with induced maps \( \imath_i:\mm_i\to\mm_{i+1} \), \( *_{\mm_{i+1}}:M_{i+1}\to\mm_{i+1} \) and proof of commutativity \( h_{i+1} \).
\end{enumerate}
A \defemph{cellular type} \( \mm \) is a sequence of types \( \mm_0\xrightarrow[]{\imath_0}\mm_1\xrightarrow[]{\imath_1}\cdots\xrightarrow[]{\imath_{n-1}}\mm_n \), together with a proof of existence of some simplicial complex \( M \) and a realization inducing this sequence.
\end{mydef}

\subsection{Polygons}
\label{sec:polygons}

The 1-type \( \bdsimplexn{2} \) has three vertices. In order to define the link of an arbitrary realization, we will need to have \( n \)-gons for \( n\geq 3 \). For example in Figure~\ref{fig:link} the link is a 6-gon. And since \( S^1 \) could be called a 1-gon, we will also define a 1-gon, and for completeness a 2-gon.

\begin{mydef}
Define \( C(n) \), \( n\geq 3 \) to be a simplicial complex with \( C(n)_{0}=\{v_1, \ldots, v_n\} \) and edges \[ C(n)_1=\{e_1=\{v_1,v_2\}, \ldots, e_{n-1}=\{v_{n-1}, v_n\}, e_n=\{v_n, v_0\}\}. \] We call \( C(n) \) a \defemph{polygon} or \defemph{\( n \)-gon}. The realization of \( C(n) \) will be denoted \( \ccc(n) \). When it is convenient we may refer to an \( n \)-gon by \( \gr{v_1\cdots v_n} \) and to its realization by \( \hgr{v_1\cdots v_n} \).
\end{mydef}

We also have two special polygons:
\begin{mydef}
The higher inductive type \( \so \), also denoted \( \ccc(1) \), has constructors:
\begin{align*}
\so &:\Type \\
\mathsf{base}&:\so \\
\mathsf{loop}&:\mathsf{base}=\mathsf{base}
\end{align*}
\end{mydef}

\begin{mydef}
The higher inductive type \(\ccc(2) \) has constructors:
\begin{align*}
\ccc(2) &:\Type \\
v_1, v_2&:\ccc(2) \\
\ell_{12}, r_{21}&:v_1=v_2\\
\end{align*}
\end{mydef}

\begin{mylemma}\label{lem:addpoints}
\( \ccc(2)\simeq \ccc(1) \) and in fact \( \ccc(n)\simeq \ccc(n-1) \).
\end{mylemma}
\begin{myproofnonqed}
(Compare to \cite{hottbook} Lemma 6.5.1.)
\end{myproofnonqed}
First we will define \( f:\ccc(2)\to \ccc(1) \) and \( g:\ccc(1)\to \ccc(2) \), then prove they are inverses.
\begin{align*}
f(v_1)=f(v_2)&=\base &\quad g(\base)&=v_1\\
f(\ell_{12})&=\loopo&\quad g(\loopo)&=\ell_{12}\cdot r_{21}\\
f(r_{21}) &= \refl_{\base}& & \\
\end{align*}

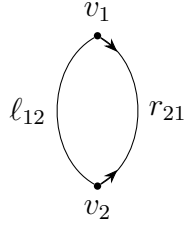
\begin{figure}[h]
\centering
\begin{tikzpicture}
\tikzset{arrow/.style={-{Stealth[scale=1.1]}}}
\tikzset{oo/.style={circle, scale=0.25, fill=black}}
\tikzset{ooo/.style={circle, scale=0.25, fill=none}}
\node[oo, label=above:\( v_1 \)] (V1) at (0, 2) {};
\node[oo, label=below:\( v_2 \)] (V2) at (0, 0) {};
\node[ooo] (V21) at (0.3, 1.75) {};
\node[ooo] (V22) at (0.3, 0.25) {};
\draw (V1) edge[bend right=60, swap, "\( \ell_{12} \)"] (V2);
\draw (V1) edge[bend left=60, "\( r_{21} \)"] (V2);
\draw[arrow] (V1) edge[bend left=5] (V21);
\draw[arrow] (V2) edge[bend right=5] (V22);
\end{tikzpicture}
\caption{We imagine shrinking \( r_{21} \) down to become \( \refl_\base \) in \( S^1 \).}
\label{fig:shrink}
\end{figure}

We need to show that \( f\circ g\sim \id_{\ccc(1)} \) and \( g\circ f\sim\id_{\ccc(2)} \).
Think of \( f \) as sliding \( v_2 \) and \( v_1 \) towards each other along \( r_{21} \) coalescing into just \( v_1 \), as in Figure~\ref{fig:shrink}. This may help understand why the somewhat intricate proof is working.

We need terms \( p:\pit{a:\ccc(1)}f(g(a))=a \) and \( q:\pit{a:\ccc(2)}g(f(a))=a \). We will proceed by induction, defining appropriate paths on point constructors and then checking a condition on path constructors that confirms that the built-in transport of these type families respects the definition on points.

Looking first at \( g\circ f \), which shrinks \( r_{21} \), we have the following data to work with:
\begin{align*}
g(f(v_1))=g(f(v_2))&=v_1\\
g(f(\ell_{12}))&=\ell_{12}\cdot r_{21}\\
g(f(r_{21})) &= \refl_{v_1}.
\end{align*}
We then need to supply a homotopy from this data to \( \id_{\ccc(2)} \), which consists of a section and pathovers over \( \ccc(2) \):
\begin{align*}
p_1&:g(f(v_1))=v_1\\
p_2&:g(f(v_1))=v_2\\
H_\ell&:\tr(\ell_{12})(p_1)=p_2\\
H_r&:\tr(r_{21})(p_2)=p_1.
\end{align*}
which simplifies to
\begin{align*}
p_1&:v_1=v_1\\
p_2&:v_1=v_2\\
H_\ell&:g(f(\ell_{12}))^{-1}\cdot p_1\cdot \ell_{12}=p_2\\
H_r&:=g(f(r_{21}))^{-1}\cdot p_2\cdot r_{21}= p_1
\end{align*}
and then to 
\begin{align*}
p_1&:v_1=v_1\\
p_2&:v_1=v_2\\
H_\ell&:(\ell_{12}\cdot r_{21})^{-1}\cdot p_1\cdot \ell_{12}=p_2\\
H_r&:\refl_{v_1}\cdot p_2\cdot r_{21}= p_1
\end{align*}

To solve all of these constraints we can choose \( p_1\defeq\refl_{v_1} \), which by consulting either \( H_\ell \) or \( H_r \) requires that we take \( p_2\defeq{r_{21}}^{-1}\).

Now examining \( f\circ g \), we have
\begin{align*}
f(g(\base))&=\base&\\
f(g(\loopo))&=f(\ell_{12}\cdot r_{21})=\loopo
\end{align*}
and so we have an easy proof that this is the identity.

The proof of the more general case \( \ccc(n) \simeq \ccc(n-1)\) is very similar. Take the maps \( f:\ccc(n)\to \ccc(n-1) \), \( g:\ccc(n-1)\to \ccc(n) \) to be
\begin{align*}
f(v_i)=v_i&\quad(i=1,\ldots,n-1) & g(v_i)&=v_i&\quad(i=1,\ldots,n-1)\\
f(v_n)=v_1&\quad& g(e_{i,i+1})&=e_{i,i+1}&\quad(i=1,\ldots,n-2)\\
f(e_{i,i+1})=e_{i,i+1}&\quad(i=1,\ldots,n-1)& g(e_{n-1,1})&=e_{n-1,n}\cdot e_{n,1}&\\
f(e_{n-1,n})=e_{n-1,1}&&&&\\
f(e_{n,1})=\refl_{v_1}&&&&
\end{align*}
where \( f \) should be thought of as shrinking \( e_{n,1} \) so that \( v_n \) coalesces into \( v_1 \).

The proof that \( g\circ f\sim\id_{\ccc(n)} \) proceeds as follows: the composition is definitionally the identity except 
\begin{align*}
g(f(v_n))&=v_1\\
g(f(e_{n-1,n}))&=e_{n-1,n}\cdot e_{n,1}\\
g(f(e_{n,1}))&= \refl_{v_1}.
\end{align*}
Guided by our previous experience we choose \( {e_{n,1}}^{-1}:g(f(v_n))=v_n \), and define the pathovers by transport.

The proof that \( f\circ g\sim\id_{\ccc(n-1)} \) requires only noting that \( f(g(e_{n-1,1}))=f(e_{n-1,n}\cdot e_{n,1})=e_{n-1,1}\cdot\refl_{v_1}=e_{n-1,1} \).\qed

\begin{mycor}
\label{cor:gon}
All polygons are equivalent to \( \so \), i.e. we have terms \( e_n:\ccc(n)=S^1 \), and hence we have constructed a map from the unit type \( (\ccc(n), ||e_n||_{-1}):\unit\to \EMzo \).
\end{mycor}
\begin{myproof}
The proofs in Lemma~\ref{lem:addpoints} can be concatenated to give \( \ccc(n)\to\ccc(n-1)\to\cdots\to\ccc(2)\to S^1 \).
\end{myproof}

\begin{mydef}
\label{def:rotation}
Let \( R:\gr{v_1\cdots v_n}\to \gr{v_1\cdots v_n} \) (for ``rotation'') be the map sending \( v_i\mapsto v_{i+1} \) and \( v_n\mapsto v_0 \). This map clearly sends edges to edges, and so is a map of simplicial complexes, and extends to a map \( \hgr{R}:\hgr{v_1\cdots v_n}\to \hgr{v_1\cdots v_n} \).
\end{mydef}

The homotopical realization \( \hgr{R} \) has a path to the identity:
\begin{mylemma}
\label{lem:rotation}
The map \( \hgr{R}: \hgr{v_1\cdots v_n}\to \hgr{v_1\cdots v_n}\) is connected to \( \refl_{\hgr{v_1\cdots v_n}} \) by a homotopy \( H_R:\pit{x:\hgr{v_1\cdots v_n}}x=\hgr{R}(x) \).
\end{mylemma}
\begin{myproof}
If \( x \) is a vertex, take \( H_R(x) \) to be the obvious unique edge back to the starting vertex. This extends in the obvious functorial way to edges.
\end{myproof}

We wrote the homotopy \( H_R(x) \) as starting at \( x \) because it feels like a record of a time-based process of applying \( R \). We will rely on this convention when we define flatness.

\subsection{Surfaces}
We will eventually focus on 2-dimensional simplicial complexes in this note, and our running example which begins in the next section is 2-dimensional. We have a simple way to define an orientation in this dimension, which we provide here. We will touch on the relationship between this classical definition and the definition of our HoTT classifying space when we discuss vector fields in Section~\ref{sec:vector_field_def}.

\begin{mydef}
An \defemph{orientation} \( \mathscr{O} \) of a 2-dimensional simplicial complex \( M=[M_0, M_1, M_2] \) is an equivalence class of ordering maps \( V:M_2\to \mathsf{List} \) to the type of ordered lists, where \( V \) orders the vertices of every 2-face. If \( F:M_2 \) is a face and \( F=\{a, b, c\} \) are its vertices, then \( V(F) \) is a list of all the vertices, e.g. \( [b, c, a] \). Two orderings \( V, V' \) are said to have the \defemph{same orientation} if it is true for every face \( F \) that \( V(F) \) and \( V'(F) \) differ by a cyclic permutation. The inclusion \( e\subset F \) of an edge in a face induces an ordering of the vertices of the edge, called \defemph{the induced order of \( e\subset F \)}.
\end{mydef}

\subsection{The octahedron model of the sphere}
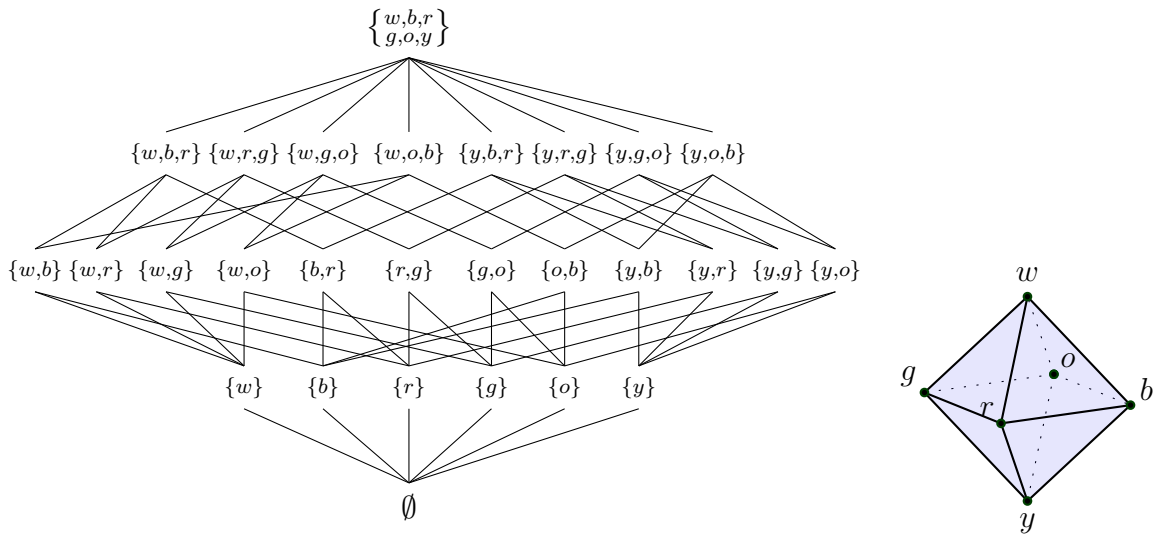
\begin{figure}[h]
\centering
\begin{tikzpicture}
    \matrix (A) [matrix of math nodes, row sep=1cm, column sep=-.2cm]
    { 
       ~ &  ~ & ~ & ~ & ~ & \left\{\substack{{w, b, r}\\ {g, o, y}}\right\} \\  
  ~ &  ~ & \scriptstyle\{w, b, r\} & \scriptstyle\{w, r, g\}  & \scriptstyle\{w, g, o\} & \scriptstyle\{w, o, b\} & \scriptstyle\{y, b, r\} & \scriptstyle\{y, r, g\}  & \scriptstyle\{y, g, o\} & \scriptstyle\{y, o, b\}\\
  \scriptstyle\{w, b\} & \scriptstyle\{w,r\}  & \scriptstyle\{w,g\} & \scriptstyle\{w,o\} & \scriptstyle\{b, r\} & \scriptstyle\{r, g\}  & \scriptstyle\{g, o\} & \scriptstyle\{o, b\} & \scriptstyle\{y, b\} & \scriptstyle\{y,r\}  & \scriptstyle\{y,g\} & \scriptstyle\{y,o\}\\
  ~ & ~ & ~ &  \scriptstyle\{w\}  & \scriptstyle\{b\} & \scriptstyle\{r\} & \scriptstyle\{g\} & \scriptstyle\{o\} & \scriptstyle\{y\}\\
      ~ & ~ &  ~ & ~ & ~ & \emptyset \\
    };
    \draw (A-1-6.south)--(A-2-3.north);
    \draw (A-1-6.south)--(A-2-4.north);
    \draw (A-1-6.south)--(A-2-5.north);
    \draw (A-1-6.south)--(A-2-6.north);
    \draw (A-1-6.south)--(A-2-7.north);
    \draw (A-1-6.south)--(A-2-8.north);
    \draw (A-1-6.south)--(A-2-9.north);
    \draw (A-1-6.south)--(A-2-10.north);

    \draw (A-2-3.south)--(A-3-1.north);
    \draw (A-2-3.south)--(A-3-2.north);
    \draw (A-2-3.south)--(A-3-5.north);

    \draw (A-2-4.south)--(A-3-2.north);
    \draw (A-2-4.south)--(A-3-3.north);
    \draw (A-2-4.south)--(A-3-6.north);

    \draw (A-2-5.south)--(A-3-3.north);
    \draw (A-2-5.south)--(A-3-4.north);
    \draw (A-2-5.south)--(A-3-7.north);

    \draw (A-2-6.south)--(A-3-4.north);
    \draw (A-2-6.south)--(A-3-1.north);
    \draw (A-2-6.south)--(A-3-8.north);

    \draw (A-2-7.south)--(A-3-5.north);
    \draw (A-2-7.south)--(A-3-10.north);
    \draw (A-2-7.south)--(A-3-9.north);

    \draw (A-2-8.south)--(A-3-6.north);
    \draw (A-2-8.south)--(A-3-11.north);
    \draw (A-2-8.south)--(A-3-10.north);

    \draw (A-2-9.south)--(A-3-7.north);
    \draw (A-2-9.south)--(A-3-12.north);
    \draw (A-2-9.south)--(A-3-11.north);

    \draw (A-2-10.south)--(A-3-8.north);
    \draw (A-2-10.south)--(A-3-9.north);
    \draw (A-2-10.south)--(A-3-12.north);

    \draw (A-3-1.south)--(A-4-4.north);
    \draw (A-3-1.south)--(A-4-5.north);

    \draw (A-3-2.south)--(A-4-4.north);
    \draw (A-3-2.south)--(A-4-6.north);

    \draw (A-3-3.south)--(A-4-4.north);
    \draw (A-3-3.south)--(A-4-7.north);

    \draw (A-3-4.south)--(A-4-4.north);
    \draw (A-3-4.south)--(A-4-8.north);

    \draw (A-3-5.south)--(A-4-5.north);
    \draw (A-3-5.south)--(A-4-6.north);

    \draw (A-3-6.south)--(A-4-6.north);
    \draw (A-3-6.south)--(A-4-7.north);

    \draw (A-3-7.south)--(A-4-7.north);
    \draw (A-3-7.south)--(A-4-8.north);

    \draw (A-3-8.south)--(A-4-8.north);
    \draw (A-3-8.south)--(A-4-5.north);

    \draw (A-3-9.south)--(A-4-9.north);
    \draw (A-3-9.south)--(A-4-5.north);

    \draw (A-3-10.south)--(A-4-9.north);
    \draw (A-3-10.south)--(A-4-6.north);

    \draw (A-3-11.south)--(A-4-9.north);
    \draw (A-3-11.south)--(A-4-7.north);

    \draw (A-3-12.south)--(A-4-9.north);
    \draw (A-3-12.south)--(A-4-8.north);

    \draw (A-4-4.south)--(A-5-6.north);
    \draw (A-4-5.south)--(A-5-6.north);
    \draw (A-4-6.south)--(A-5-6.north);
    \draw (A-4-7.south)--(A-5-6.north);
    \draw (A-4-8.south)--(A-5-6.north);
    \draw (A-4-9.south)--(A-5-6.north);
\end{tikzpicture}
\begin{tikzpicture}%
  [x={(-0.860769cm, -0.121512cm)},
  y={(0.508996cm, -0.205391cm)},
  z={(-0.000053cm, 0.971107cm)},
  scale=1,
  back/.style={loosely dotted, thin},
  edge/.style={black, thick},
  facet/.style={fill=blue!95!black,fill opacity=0.1},
  vertex/.style={inner sep=1pt,circle,draw=green!25!black,fill=black,thick}]
\coordinate (-1, -1, 0) at (-1, -1, 0);
\coordinate (-1, 1, 0) at (-1, 1, 0);
\coordinate (0, 0, -1) at (0, 0, -1);
\coordinate (0, 0, 1) at (0, 0, 1);
\coordinate (1, -1, 0) at (1, -1, 0);
\coordinate (1, 1, 0) at (1, 1, 0);
\draw[edge,back] (-1, -1, 0) -- (-1, 1, 0);
\draw[edge,back] (-1, -1, 0) -- (0, 0, -1.4);
\draw[edge,back] (-1, -1, 0) -- (0, 0, 1.4);
\draw[edge,back] (-1, -1, 0) -- (1, -1, 0);
\node[vertex] at (-1, -1, 0)     {};
\fill[facet] (1, 1, 0) -- (0, 0, -1.4) -- (1, -1, 0) -- cycle {};
\fill[facet] (1, 1, 0) -- (0, 0, 1.4) -- (1, -1, 0) -- cycle {};
\fill[facet] (1, 1, 0) -- (-1, 1, 0) -- (0, 0, 1.4) -- cycle {};
\fill[facet] (1, 1, 0) -- (-1, 1, 0) -- (0, 0, -1.4) -- cycle {};
\draw[edge] (-1, 1, 0) -- (0, 0, -1.4);
\draw[edge] (-1, 1, 0) -- (0, 0, 1.4);
\draw[edge] (-1, 1, 0) -- (1, 1, 0);
\draw[edge] (0, 0, -1.4) -- (1, -1, 0);
\draw[edge] (0, 0, -1.4) -- (1, 1, 0);
\draw[edge] (0, 0, 1.4) -- (1, -1, 0);
\draw[edge] (0, 0, 1.4) -- (1, 1, 0);
\draw[edge] (1, -1, 0) -- (1, 1, 0);
\begin{scope}[nodes=vertex]
\node[label=above right:\( b \)] at (-1, 1, 0)     {};
\node[label=below:\( y \)] at (0, 0, -1.4)     {};
\node[label=above:\( w \)] at (0, 0, 1.4)     {};
\node[label=above left:\( g \)] at (1, -1, 0)     {};
\node[label=above left:\( r \)] at (1, 1, 0)     {};
\node[label=above right:\( o \)] at (-1, -1, 0)     {};
\end{scope}
\end{tikzpicture}
\caption{The Hasse diagram of the simplicial complex \( O \), and a possible realization. The row of singletons in the Hasse diagram is \( O_0 \) and above it are \( O_1 \) and \( O_2 \).}
\label{fig:hasse_octohedron}
\end{figure}

We will create our first combinatorial surface, an octahedron. In \( \simcomp \) the combinatorial data of the faces can be represented with a \emph{Hasse diagram}, which shows the poset of inclusions in a graded manner, with a special top and bottom element. The top element is part of the diagram only, and not part of the simplicial complex, else it would provide a filling cell. We give an octahedron \( O=[O_0, O_1, O_2] \) in Figure~\ref{fig:hasse_octohedron}. The names of the vertices are short for white, yellow, blue, red, green, and orange, the colors of the faces of a Rubik's cube. The octahedron is the dual of the cube, with each vertex corresponding to a face.

We can realize \( O=[O_0, O_1, O_2] \) as a cellular type \( \oo \).

\begin{mylemma}
\label{lem:octahedron_sphere}
There is an equivalence \( \oo_2\simeq S^2 \).
\end{mylemma}
\begin{myproof}Omitted.\end{myproof}

\clearpage
\section{Bundles, connections, and curvature}
\label{sec:connections}
Bundles are simply maps into the universe. By using the extra cellular structure and the even more detailed combinatorial structure of realizations, we can identify inside of HoTT some additional classical definitions.

\subsection{Definitions}
Having the cellular structure allows us to define connections.
\begin{mydef}
\label{def:connection}
If \( \mm\defeq \mm_0\xrightarrow[]{\imath_0}\cdots\xrightarrow[]{\imath_{n-1}}\mm_n \) is a cellular type and \( f_k:\mm_k\to\uni \) are type families on each skeleton such that all the triangles commute in the diagram:
\end{mydef}
\begin{center}
\begin{tikzcd}
  {\mm_0} & {\mm_1} & {\mm_2} & \cdots & {\mm_n} \\
  && {\mathcal{U}}
  \arrow["{\imath_0}", from=1-1, to=1-2]
  \arrow["{f_0}", from=1-1, to=2-3]
  \arrow["{\imath_1}", from=1-2, to=1-3]
  \arrow["{f_1}", from=1-2, to=2-3]
  \arrow["{\imath_2}", from=1-3, to=1-4]
  \arrow["{f_2}", from=1-3, to=2-3]
  \arrow["{\imath_{n-1}}", from=1-4, to=1-5]
  \arrow["f_n"', from=1-5, to=2-3]
\end{tikzcd}
\end{center}
then we say
\begin{itemize}
\item The map \( f_k \) is a \defemph{\( k \)-bundle} on \( \mm \).
\item The pair given by the map \( f_k \) and the proof \( f_k\circ \imath_{k-1}=f_{k-1} \) that \( f_k \) extends \( f_{k-1} \) is called a \defemph{\( k \)-connection on the \( (k-1) \)-bundle \( f_{k-1} \)}.
\end{itemize}

Having the additional structure of a simplicial complex allows us to define curvature, which is a local concept.

\begin{mydef}
If \( \mm \) is the realization of a simplicial complex, such that for each pushout defining \( \mm_k \) we have the diagram
\end{mydef}
\begin{center}
\begin{tikzcd}
  {M_k\times \bdsimplexn{k}} & {M_k} \\
  {\mathbb{M}_{k-1}} & {\mathbb{M}_k} \\
  & {\mathcal{U}}
  \arrow["{\mathrm{pr}_1}", from=1-1, to=1-2]
  \arrow["{\mathbb{A}_{k-1}}"', from=1-1, to=2-1]
  \arrow["{*_{\mathbb{M}_k}}", from=1-2, to=2-2]
  \arrow["{h_k}", shorten <=10pt, shorten >=10pt, Rightarrow, from=2-1, to=1-2]
  \arrow["{\imath_{k-1}}", from=2-1, to=2-2]
  \arrow[""{name=0, anchor=center, inner sep=0}, "{f_{k-1}}"', from=2-1, to=3-2]
  \arrow["\ulcorner"{anchor=center, pos=0.125, rotate=180}, draw=none, from=2-2, to=1-1]
  \arrow["{f_k}", from=2-2, to=3-2]
  \arrow[shorten >=3pt, Rightarrow, from=2-2, to=0]
\end{tikzcd}
\end{center}

the outer square of which restricts on each face to the diagram

\begin{center}
\begin{tikzcd}
  {\{F\}\times \bdsimplexn{k}} & \{F\} \\
  {\mathbb{M}_{k-1}} & {\mathcal{U}}
  \arrow["{\pr_1}", from=1-1, to=1-2]
  \arrow["{\mathbb{A}_{k-1}}"', from=1-1, to=2-1]
  \arrow["{*_{\mathbb{M}_k}}", from=1-2, to=2-2]
  \arrow["{\flat_k}", shorten <=11pt, shorten >=11pt, Rightarrow, from=1-2, to=2-1]
  \arrow[from=2-1, to=2-2]
\end{tikzcd}
\end{center}
then we say the filler \( \flat_k \) is called a \defemph{flatness structure for the face \( m_k \)}, and its ending path is called \defemph{curvature at the face \( F \)}.

Although all the pushout diagrams are introducing hub points for each edge and face, in practice we will ignore those. But in the case of flatness, which we will use extensively, we should clarify what its type is when we are ignoring hubs. If we have a face \( F \) with vertices \( v_1, v_2, v_3 \) and paths \( e_{ij}:v_i=v_j \) that are realizations of edges, and the loop \( \ell_F\defeq e_{12}\cdot e_{23}\cdot e_{31}\) then we will write \(\flat(\ell_F):\id_{f_1(v_1)}=f_1(\ell_F)\) for a 2-path filling this loop.

The definitions can be digested to give
\begin{mylemma}
\label{lem:extend_faces}
Given \( f_{k-1} \) as above, a \( k \)-connection exists if and only if there exists a flatness structure for each \( k \)-face.
\end{mylemma}

On a 2-dimensional cellular type \( \mm\defeq\mm_0\to\mm_1\to\mm_2 \) the terminology works out as follows: a 0-bundle on \( \mm \) is a map \( T_0:\mm_0\to\uni \). A 1-connection on \( T_0 \) is an extension \( T_1:\mm_1\to\uni \). A 2-connection on \( T_1 \) is an extension \( T_2:\mm_2\to\uni \). Classically, a 1-connection that extends to a 2-connection is called \emph{flat}.

\subsection{Flat connections as local trivializations}
\label{sec:localtriv}
This section can be viewed as an extended remark. The observation we want to make is that the data of a 2-bundle on a realization of a 2-dimensional simplicial complex is related to the construction of local trivializations: the fiber at one vertex can be extended throughout a single face coherently, using the connection (the extension of the classifying map to the edges) to specify isomorphisms with the fibers at the other points, and the higher connections to establish commutativity between these. 

We introduce a notation more suitable for the algebra of charts and overlaps: denote the fiber at \( v_i \) by \( T_i \) and denote transport along \( e_{ij}:v_i=v_j \) by \( T_{ji}:T_i\to T_j \). The indices are ordered from right to left, which is compatible with function composition notation. Denote the inverse function by swapping indices: \( T_{ij}\defeq T_{ji}^{-1} \). Assume we have some fixed isomorphism \( T_i=S^1 \), and to avoid composing everything with this function we will assume it is \( \id \). In the diagram below we see the data arranged so that our bundles fibers are on the left, and the fiber of a trivial bundle is on the right.
\begin{center}
\begin{tikzcd}
  {T_i} & {S^1} \\
  {T_j} & {S^1} \\
  {T_k} & {S^1}
  \arrow["{\mathrm{id}}", equals, from=1-1, to=1-2]
  \arrow["{T_{ji}}"', from=1-1, to=2-1]
  \arrow[""{name=0, anchor=center, inner sep=0}, "{T_{ki}}"', curve={height=40pt}, from=1-1, to=3-1]
  \arrow[equals, from=1-2, to=2-2]
  \arrow[curve={height=-40pt}, equals, from=1-2, to=3-2]
  \arrow["{T_{ij}}"', from=2-1, to=2-2]
  \arrow["{T_{kj}}"', from=2-1, to=3-1]
  \arrow[equals, from=2-2, to=3-2]
  \arrow["{T_{ij}T_{jk}}"', from=3-1, to=3-2]
  \arrow["{\flat_{ijk}}"', shorten >=4pt, Rightarrow, from=2-1, to=0]
\end{tikzcd}
\end{center}
The two middle squares commute definitionally. Call these two squares together the back face. The left triangle is filled by the flatness structure on the face, and the right triangular filler is trivial. There is also a filler needed for the front, i.e. the outer square. This requires proving that \( T_{ki}T_{ij}T_{jk}=\id \), which is supplied by the flatness structure. There is also a 3-cell filling the interior of this prism, mapping the back face plus the left triangle filler to the front face plus the right triangle filler. These two faces both consist of one or two identities concatenated with the flatness structure, and so the 3-cell is definitional.

This relationship between flatness structure and local triviality of a chart can be compared to the classical result that on a paracompact, simply connected manifold (such as a single chart), a connection on a principal bundle is flat if and only if the bundle is trivial. See for example \cite{kobayashinomizu} Corollary 9.2.

\subsection{The tangent bundle of the sphere}
We will build up a map \( T \) out of \( \oo_0\to\oo_1\to\oo_2 \) which is meant to be a model of the tangent bundle of the sphere. The link function will serve as our approximation to the tangent space. Taking the link of a vertex gives us a map from vertices to polygons, so the codomain is \( \EMzo \).

If \( \{b, r, g, o\} \) are four vertices in \( \oo \), the notation \( \hgr{brgo} \) refers to the 4-gon spanned by these four vertices and the edges of \( \oo \) that connect them to each other.

\begin{mydef}
\( T_0\defeq\link:\oo_0\to\EMzo \) is given by:
\begin{align*}
\link(w) &= \hgr{brgo} & \link(r) &= \hgr{wbyg} \\
\link(y) &= \hgr{bogr} & \link(g) &= \hgr{wryo} \\
\link(b) &= \hgr{woyr} & \link(o) &= \hgr{wgyb}
\end{align*}

We chose these orderings for the vertices in the link, by visualizing standing at the given vertex as if it were the north pole, then looking south and enumerating the link in clockwise order, starting from \( w \) if possible, else \( b \).
\end{mydef}

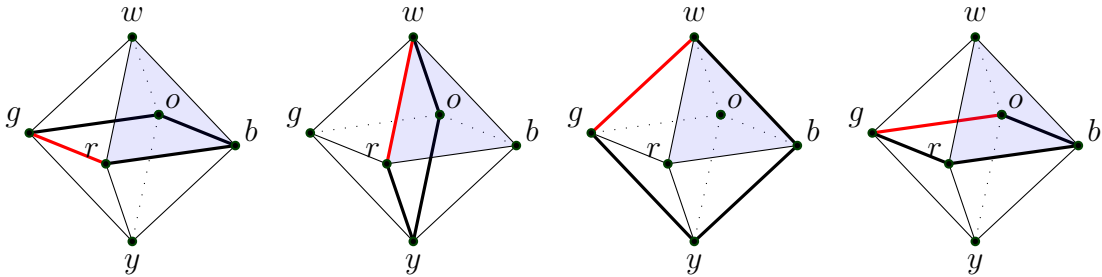
\begin{figure}[h]
\centering
\begin{tikzpicture}%
  [x={(-0.860769cm, -0.121512cm)},
  y={(0.508996cm, -0.205391cm)},
  z={(-0.000053cm, 0.971107cm)},
  scale=1,
  eqback/.style={very thick},
  back/.style={loosely dotted, thin},
  eqedge/.style={very thick},
  edge/.style={black, thin},
  r/.style={red},
  facet/.style={fill=blue!95!black,fill opacity=0.1},
  vertex/.style={inner sep=1pt,circle,draw=green!25!black,fill=black,thick}]
\coordinate (-1, -1, 0) at (-1, -1, 0);
\coordinate (-1, 1, 0) at (-1, 1, 0);
\coordinate (0, 0, -1) at (0, 0, -1);
\coordinate (0, 0, 1) at (0, 0, 1);
\coordinate (1, -1, 0) at (1, -1, 0);
\coordinate (1, 1, 0) at (1, 1, 0);
\draw[edge,eqback] (-1, -1, 0) -- (-1, 1, 0);
\draw[edge,back] (-1, -1, 0) -- (0, 0, -1.4);
\draw[edge,back] (-1, -1, 0) -- (0, 0, 1.4);
\draw[edge,eqback] (1, -1, 0) -- (-1, -1, 0);
\node[vertex] at (-1, -1, 0)     {};
\fill[facet] (1, 1, 0) -- (-1, 1, 0) -- (0, 0, 1.4) -- cycle {};
\draw[edge] (-1, 1, 0) -- (0, 0, -1.4);
\draw[edge] (-1, 1, 0) -- (0, 0, 1.4);
\draw[eqedge] (-1, 1, 0) -- (1, 1, 0);
\draw[edge] (0, 0, -1.4) -- (1, -1, 0);
\draw[edge] (0, 0, -1.4) -- (1, 1, 0);
\draw[edge] (0, 0, 1.4) -- (1, -1, 0);
\draw[edge] (0, 0, 1.4) -- (1, 1, 0);
\draw[r,eqedge] (1, 1, 0) -- (1, -1, 0);
\begin{scope}[nodes=vertex]
\node[label=above right:\( b \)] at (-1, 1, 0)     {};
\node[label=below:\( y \)] at (0, 0, -1.4)     {};
\node[label=above:\( w \)] at (0, 0, 1.4)     {};
\node[label=above left:\( g \)] at (1, -1, 0)     {};
\node[label=above left:\( r \)] at (1, 1, 0)     {};
\node[label=above right:\( o \)] at (-1, -1, 0)     {};
\end{scope}
\end{tikzpicture}
\begin{tikzpicture}%
  [x={(-0.860769cm, -0.121512cm)},
  y={(0.508996cm, -0.205391cm)},
  z={(-0.000053cm, 0.971107cm)},
  scale=1,
  eqback/.style={very thick},
  back/.style={loosely dotted, thin},
  eqedge/.style={very thick},
  r/.style={red},
  edge/.style={black, thin},
  facet/.style={fill=blue!95!black,fill opacity=0.1},
  vertex/.style={inner sep=1pt,circle,draw=green!25!black,fill=black,thick}]
\coordinate (-1, -1, 0) at (-1, -1, 0);
\coordinate (-1, 1, 0) at (-1, 1, 0);
\coordinate (0, 0, -1) at (0, 0, -1);
\coordinate (0, 0, 1) at (0, 0, 1);
\coordinate (1, -1, 0) at (1, -1, 0);
\coordinate (1, 1, 0) at (1, 1, 0);
\draw[edge,back] (-1, -1, 0) -- (-1, 1, 0);
\draw[edge,eqback] (-1, -1, 0) -- (0, 0, -1.4);
\draw[edge,eqback] (0, 0, 1.4) -- (-1, -1, 0);
\draw[edge,back] (1, -1, 0) -- (-1, -1, 0);
\node[vertex] at (-1, -1, 0)     {};
\fill[facet] (1, 1, 0) -- (-1, 1, 0) -- (0, 0, 1.4) -- cycle {};
\draw[edge] (-1, 1, 0) -- (0, 0, -1.4);
\draw[edge] (-1, 1, 0) -- (0, 0, 1.4);
\draw[edge] (-1, 1, 0) -- (1, 1, 0);
\draw[edge] (0, 0, -1.4) -- (1, -1, 0);
\draw[eqedge] (0, 0, -1.4) -- (1, 1, 0);
\draw[edge] (0, 0, 1.4) -- (1, -1, 0);
\draw[r,eqedge] (1, 1, 0) -- (0, 0, 1.4) ;
\draw[edge] (1, 1, 0) -- (1, -1, 0);
\begin{scope}[nodes=vertex]
\node[label=above right:\( b \)] at (-1, 1, 0)     {};
\node[label=below:\( y \)] at (0, 0, -1.4)     {};
\node[label=above:\( w \)] at (0, 0, 1.4)     {};
\node[label=above left:\( g \)] at (1, -1, 0)     {};
\node[label=above left:\( r \)] at (1, 1, 0)     {};
\node[label=above right:\( o \)] at (-1, -1, 0)     {};
\end{scope}
\end{tikzpicture}
\begin{tikzpicture}%
  [x={(-0.860769cm, -0.121512cm)},
  y={(0.508996cm, -0.205391cm)},
  z={(-0.000053cm, 0.971107cm)},
  scale=1,
  eqback/.style={very thick},
  back/.style={loosely dotted, thin},
  eqedge/.style={very thick},
  r/.style={red},
  edge/.style={black, thin},
  facet/.style={fill=blue!95!black,fill opacity=0.1},
  vertex/.style={inner sep=1pt,circle,draw=green!25!black,fill=black,thick}]
\coordinate (-1, -1, 0) at (-1, -1, 0);
\coordinate (-1, 1, 0) at (-1, 1, 0);
\coordinate (0, 0, -1) at (0, 0, -1);
\coordinate (0, 0, 1) at (0, 0, 1);
\coordinate (1, -1, 0) at (1, -1, 0);
\coordinate (1, 1, 0) at (1, 1, 0);
\draw[edge,back] (-1, -1, 0) -- (-1, 1, 0);
\draw[edge,back] (-1, -1, 0) -- (0, 0, -1.4);
\draw[edge,back] (-1, -1, 0) -- (0, 0, 1.4);
\draw[edge,back] (1, -1, 0) -- (-1, -1, 0);
\node[vertex] at (-1, -1, 0)     {};
\fill[facet] (1, 1, 0) -- (-1, 1, 0) -- (0, 0, 1.4) -- cycle {};
\draw[eqedge] (-1, 1, 0) -- (0, 0, -1.4);
\draw[eqedge] (0, 0, 1.4) -- (-1, 1, 0);
\draw[edge] (-1, 1, 0) -- (1, 1, 0);
\draw[eqedge] (0, 0, -1.4) -- (1, -1, 0);
\draw[edge] (0, 0, -1.4) -- (1, 1, 0);
\draw[r,eqedge] (1, -1, 0) -- (0, 0, 1.4);
\draw[edge] (0, 0, 1.4) -- (1, 1, 0);
\draw[edge] (1, 1, 0) -- (1, -1, 0);
\begin{scope}[nodes=vertex]
\node[label=above right:\( b \)] at (-1, 1, 0)     {};
\node[label=below:\( y \)] at (0, 0, -1.4)     {};
\node[label=above:\( w \)] at (0, 0, 1.4)     {};
\node[label=above left:\( g \)] at (1, -1, 0)     {};
\node[label=above left:\( r \)] at (1, 1, 0)     {};
\node[label=above right:\( o \)] at (-1, -1, 0)     {};
\end{scope}
\end{tikzpicture}
\begin{tikzpicture}%
  [x={(-0.860769cm, -0.121512cm)},
  y={(0.508996cm, -0.205391cm)},
  z={(-0.000053cm, 0.971107cm)},
  scale=1,
  eqback/.style={very thick},
  back/.style={loosely dotted, thin},
  eqedge/.style={very thick},
  edge/.style={black, thin},
  r/.style={red},
  facet/.style={fill=blue!95!black,fill opacity=0.1},
  vertex/.style={inner sep=1pt,circle,draw=green!25!black,fill=black,thick}]
\coordinate (-1, -1, 0) at (-1, -1, 0);
\coordinate (-1, 1, 0) at (-1, 1, 0);
\coordinate (0, 0, -1) at (0, 0, -1);
\coordinate (0, 0, 1) at (0, 0, 1);
\coordinate (1, -1, 0) at (1, -1, 0);
\coordinate (1, 1, 0) at (1, 1, 0);
\draw[edge,eqback] (-1, -1, 0) -- (-1, 1, 0);
\draw[edge,back] (-1, -1, 0) -- (0, 0, -1.4);
\draw[edge,back] (-1, -1, 0) -- (0, 0, 1.4);
\draw[edge,eqback,r] (1, -1, 0) -- (-1, -1, 0);
\node[vertex] at (-1, -1, 0)     {};
\fill[facet] (1, 1, 0) -- (-1, 1, 0) -- (0, 0, 1.4) -- cycle {};
\draw[edge] (-1, 1, 0) -- (0, 0, -1.4);
\draw[edge] (-1, 1, 0) -- (0, 0, 1.4);
\draw[eqedge] (-1, 1, 0) -- (1, 1, 0);
\draw[edge] (0, 0, -1.4) -- (1, -1, 0);
\draw[edge] (0, 0, -1.4) -- (1, 1, 0);
\draw[edge] (0, 0, 1.4) -- (1, -1, 0);
\draw[edge] (0, 0, 1.4) -- (1, 1, 0);
\draw[eqedge] (1, 1, 0) -- (1, -1, 0);
\begin{scope}[nodes=vertex]
\node[label=above right:\( b \)] at (-1, 1, 0)     {};
\node[label=below:\( y \)] at (0, 0, -1.4)     {};
\node[label=above:\( w \)] at (0, 0, 1.4)     {};
\node[label=above left:\( g \)] at (1, -1, 0)     {};
\node[label=above left:\( r \)] at (1, 1, 0)     {};
\node[label=above right:\( o \)] at (-1, -1, 0)     {};
\end{scope}
\end{tikzpicture}
\caption{\( \link \) for the vertices \( w, b\) and \( r \).}
\label{fig:triangle_of_equators}
\end{figure}

To extend \( T_0 \) to a function \( T_1 \) on the 1-skeleton we have some freedom. We will do something motivated by the figures we have been drawing of an octahedron embedded in 3-dimensional space. We will imagine how \( T_1 \) changes as we slide from point to point in the embedding shown in the figures. Sliding from \( w \) to \( b \) and tipping the link as we go, we see \( r\mapsto r \) and \( o\mapsto o \) because those lie on the axis of rotation. Then \( g\mapsto w \) and \( b\mapsto y \).

\begin{mydef}
Define \( T_1:\oo_1\to\EMzo \) on just the 1-skeleton by extending \( T_0 \) as follows:
Transport away from \( w \):
\begin{itemize}
\item \( T_1(wr):\hgr{brgo}\mapsto \hgr{bygw} \) (\( b, g \) fixed)
\item \( T_1(wg):\hgr{brgo}\mapsto \hgr{wryo} \)
\item \( T_1(wb):\hgr{brgo}\mapsto \hgr{yrwo} \) (\( r, o \) fixed)
\item \( T_1(wo):\hgr{brgo}\mapsto \hgr{bwgy} \)
\end{itemize}
Transport away from \( y \):
\begin{itemize}
\item \( T_1(yb):\hgr{bogr}\mapsto \hgr{woyr} \)
\item \( T_1(yr):\hgr{bogr}\mapsto \hgr{bygw} \)
\item \( T_1(yg):\hgr{bogr}\mapsto \hgr{yowr} \)
\item \( T_1(yo):\hgr{bogr}\mapsto \hgr{bwgy} \)
\end{itemize}
Transport along the equator:
\begin{itemize}
\item \( T_1(br):\hgr{woyr}\mapsto \hgr{wbyg} \) 
\item \( T_1(rg):\hgr{wbyg}\mapsto \hgr{wryo} \)
\item \( T_1(go):\hgr{wryo}\mapsto \hgr{wgyb} \)
\item \( T_1(ob):\hgr{wgyb}\mapsto \hgr{woyr} \)
\end{itemize}
\label{def:transport_oo}
\end{mydef}

It's very important to be able to visualize what \( T_1 \) does to triangular paths such as \( wb\cdot br\cdot rw \) (which circulates around the boundary of face \( wbr \)). You can see it if you imagine Figure~\ref{fig:triangle_of_equators} as the frames of a short movie. Or you can place your palm over the top of a cube and note where your fingers are pointing, then slide your hand to an equatorial face, then along the equator, then back to the top. The answer is: you come back rotated clockwise by a quarter-turn, which we saw in Definition~\ref{def:rotation} where it is called \( \hgr{R} \).

Now let's extend \( T_1 \) to all of \( \oo \) by providing values for the eight faces. The face \( wbr \) is a path from \( \refl_w \) to the concatenation \( wb\cdot br\cdot rw \), and so the image of \( wbr \) under the extended version of \( T_1 \) must be a homotopy from \( \refl_{T_1(w)} \) to \( T_1(wb\cdot br\cdot rw) \). Here \emph{there is no additional freedom}.

\begin{mydef}
\label{def:octahedron_curvature}
Define \( T_2:\oo\to\EMzo \) by extending \( T_1 \) to the faces as follows (making use of \( H_R \) from Lemma~\ref{lem:rotation}):
\begin{multicols}{2}
\begin{itemize}
\item \( T_2(wbr)=H_R \) 
\item \( T_2(wrg)=H_R \)
\item \( T_2(wgo)=H_R \)
\item \( T_2(ybo)=H_R \)
\item \( T_2(yrb)=H_R \) 
\item \( T_2(ygr)=H_R \)
\item \( T_2(yog)=H_R \)
\item \( T_2(ybo)=H_R \)
\end{itemize}
\end{multicols}
Defining these flatness structures suffices to define \( T_2 \) by Lemma~\ref{lem:extend_faces}.
\end{mydef}

\subsection{Existence of connections}
\label{sec:existence}
How confident can we be that we can always define a connection on the tangent bundle of an arbitrary simplicial complex? Two things make the octahedron example special: the link is a 4-gon at every vertex (as opposed to having a variable number of vertices), and every transport map extends to a rotation of the entire octahedron in 3-dimensional space. This imposed a coherence on the interactions of all the choices we made for the connection, which we can worry may not exist for more complex combinatorial data.

We know as a fact outside of HoTT that any combinatorial surface that has been realized as a triangulated surface embedded in 3-dimensional euclidean space can inherit the parallel transport entailed in the embedding. We could then approximate that data to arbitrary precision with enough subdivision of the fibers of \( T \).

What would a proof inside of HoTT look like? We will leave this as an open question.

\clearpage
\section{Vector fields}
\label{sec:vector_fields}
\subsection{Definition}
\label{sec:vector_field_def}
Vector fields are sections of the tangent bundle of a manifold. If \( \mm\defeq\mm_0\to\mm_1\to\mm_2 \) is a cellular type, and given an extension \( T:\mm\to\EMzo \) of the \( \link \) function, we can consider the type of sections \( \pit{x:\mm_1}T_1(x) \).

In this section and for the remainder of the note, we will assume that the bundle is principal, that is that the bundle is in fact a map \( T:\mm\to\Kzt \), the type of \( S^1 \)-torsors. This assumption is in principle not needed if \( M \) is equipped with an orientation \( \mathscr{O} \), but we are omitting a proof that such an \( \mathscr{O} \) provides a factorization of \( T \).

\begin{mydef}
Given a 2-dimensional cellular type \( \mm\defeq\mm_0\to\mm_1\to\mm_2 \) equipped with type family \( T:\mm\to\Kzt \) extending \( \link \), a \defemph{vector field} on \( \mm \) is a section on the 1-skeleton, i.e. a term \( X:\pit{x:\mm_1}T(x) \).
\end{mydef}

\begin{mynote}
The circle bundle extending \( \link \) captures the \emph{unit spheres} of the classical tangent bundle. A section of this bundle is therefore analogous to a classical \emph{nonvanishing} vector field. To allow for classical vector fields \emph{with zeros}, we are limiting the section \( X \) to the 1-skeleton of \( \mm \).
\end{mynote}

On the 0-skeleton \( X \) picks a point in each link, i.e. a neighbor of each vertex. On a path \( p:x=_\mm y \), \( X \) assigns a dependent path over \( p \), which as we know is a term \( \pi:\tr(p)(X(x))=_{T y} X(y) \). We are very interested in working with the concatenation operation on dependent paths, which we call \emph{swirling}.

\subsection{Swirling of dependent paths}
\begin{figure}[h]
\centering
\begin{tikzpicture}
  [arrow/.style={-{Stealth[scale=1.1]}}, vec/.style={ultra thick, color=black}, vectr/.style={thick, color=black}, vectrtr/.style={thick, dashed, color=black}, vectrtrtr/.style={thick, dotted, color=black}]
  \tikzset{oo/.style={circle, scale=0.6, fill=black}}
  \tikzset{ii/.style={circle, scale=0.3, fill=gray}}
  \setlength{\mylen}{3cm}
  \setlength{\mylin}{1.2cm}
    \node[oo, label=below right:\( v_1 \)] (V1) at (0, 0) {};
    \node[oo, label=below:\( v_3 \)] (V3) at (2*\mylen, 0) {};
    \node[oo, label=above:\( v_2 \)] (V2) at (\mylen, 1.732*\mylen) {};

    \draw[arrow] (V2) edge[very thick, color=teal, "\( e_{23} \)"] (V3);
    \draw[arrow] (V1) edge[very thick, color=magenta, "\( e_{12} \)"] (V2);
    \draw[arrow] (V3) edge[very thick, color=blue, "\( e_{31} \)"] (V1);
    
    \node [ii, above right=\mylin of V1] (V11) {};
    \node [ii, below right=\mylin of V1] (V14) {};
    \node [ii, below left=\mylin of  V1] (V13) {};
    \node [ii, above left= \mylin of V1] (V12) {};

    \node [left=1.3\mylin of  V1,  label=center:\( T_1 \)] {};
    \node [right=1.3\mylin of  V2,  label=center:\( T_2 \)] {};
    \node [right=1.3\mylin of  V3,  label=center:\( T_3 \)] {};

    \node [ii, above right=\mylin of V2] (V21) {};
    \node [ii, below right=\mylin of V2] (V24) {};
    \node [ii, below left=\mylin of  V2] (V23) {};
    \node [ii, above left= \mylin of V2] (V22) {};

    \node [ii, above right=\mylin of V3] (V31) {};
    \node [ii, below right=\mylin of V3] (V34) {};
    \node [ii, below left=\mylin of  V3] (V33) {};
    \node [ii, above left= \mylin of V3] (V32) {};

    \draw[dashed] (V11) -- (V12);
    \draw[dashed] (V12) -- (V13);
    \draw[dashed] (V13) -- (V14);
    \draw[dashed] (V14) -- (V11);

    \draw[dashed] (V21) -- (V22);
    \draw[dashed] (V22) -- (V23);
    \draw[dashed] (V23) -- (V24);
    \draw[dashed] (V24) -- (V21);
    
    \draw[dashed] (V31) -- (V32);
    \draw[dashed] (V32) -- (V33);
    \draw[dashed] (V33) -- (V34);
    \draw[dashed] (V34) -- (V31);
    
    \draw[arrow] (V1) edge[vec] (V11);
    \draw[arrow] (V2) edge[vectr] (V21);
    \draw[arrow] (V3) edge[vectrtr] (V34);
    \draw[arrow] (V1) edge[vectrtrtr] (V14);

    \draw[arrow] (V2) edge[vec] (V24);
    \draw[arrow] (V3) edge[vectr] (V33);
    \draw[arrow] (V1) edge[vectrtr] (V13);

    \draw[arrow] (V21) edge[thick, color=magenta] (V24);
    \draw[arrow] (V34) edge[thick, color=magenta] (V33);
    \draw[arrow] (V14) edge[thick, color=magenta] (V13);
    \draw[arrow] (V33) edge[thick, color=teal] (V32);
    \draw[arrow] (V13) edge[thick, color=teal] (V12);
    \draw[arrow] (V12) edge[thick, color=blue] (V11);

    \draw[arrow] (V3) edge[vec] (V32);
    \draw[arrow] (V1) edge[vectr] (V12);
\end{tikzpicture}
\caption{A vector field swirling clockwise around a face, in a bundle of squares. Imagine that transport along \( e_{12} \) does not rotate along the page, that transport along \( e_{23} \) rotates clockwise by 90 degrees, and that transport along \( e_{31} \) again does not rotate along the page. Thick black vectors are the vector field at a point. Thin vectors are transported once, dashed twice, and dotted three times. The vertices \( v_{ij} \) are in the tangent fibers. If you have a colorized version of the document, the colors of the arrows correspond: the red edge produced the red edge in the fibers.}
\label{fig:swirling}
\end{figure}
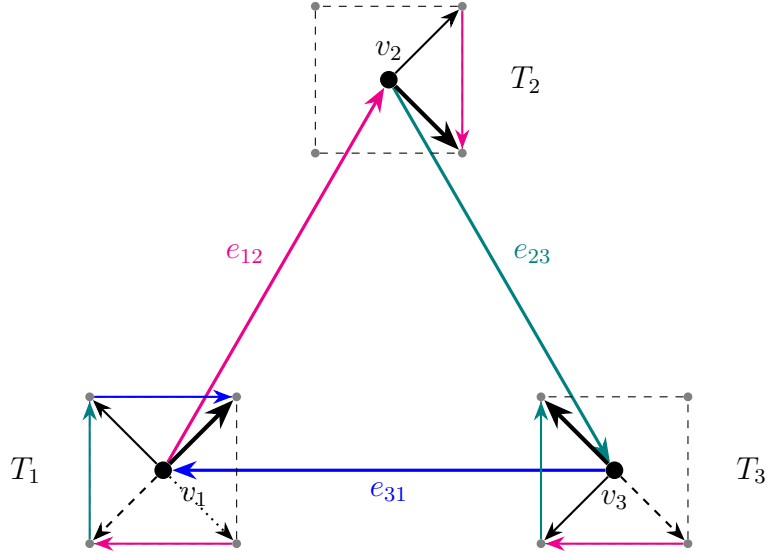

Consider the vertex \( v_1:\mm \), a face \( F \) containing vertices \( v_1, v_2, v_3 \), and the boundary path \( \ell\defeq e_{12}\cdot e_{23} \cdot e_{31} \). As we did in Section~\ref{sec:localtriv}, denote \( T(v_i) \) by \( T_i \) and \( T(e_{ij}) \) by \( T_{ji} \). The indices are swapped so that we can have expressions that respect function-composition order, such as \( T_{32}T_{21}(X_1):T_3 \). We retain the opposite convention (which we can call path-concatenation order) for the vector field, e.g. \( X_{ij} \) which is a path in \( T_j \), as these are paths. Figure~\ref{fig:vector_transport} shows in tabular form how we concatenate the dependent paths over \( e_{12}\cdot e_{23}\cdot e_{31} \). Figure~\ref{fig:swirling} shows visually a possible example.

\begin{figure}[h]
\centering
\begin{tikzcd}
  {T_1} & {T_2} & {T_3} & {T_1} \\
  &&& {T_{13}T_{32}T_{21}X_1} \\
  && {T_{32}T_{21}X_1} & {T_{13}T_{32}X_2} \\
  & {T_{21}X_1} & {T_{32}X_2} & {T_{13}X_3} \\
  {X_1} & {X_2} & {X_3} & {X_1}
  \arrow["{T_{21}}", from=1-1, to=1-2]
  \arrow["{T_{32}}", from=1-2, to=1-3]
  \arrow["{T_{13}}", from=1-3, to=1-4]
  \arrow["{T_{13}T_{32}X_{12}:}", equals, from=3-4, to=2-4]
  \arrow["{X_{12}:}"', equals, from=4-2, to=5-2]
  \arrow["{T_{32}X_{12}:}", equals, from=4-3, to=3-3]
  \arrow["{T_{13}X_{23}:}", equals, from=4-4, to=3-4]
  \arrow["{X_{23}:}", equals, from=5-3, to=4-3]
  \arrow["{X_{31}:}", equals, from=5-4, to=4-4]
\end{tikzcd}
\caption{The data in each fiber as we move around a triangle with vertices indexed 1, 2, and 3. Double lines indicate identity types between two types, and their labels are terms of this type (reading downwards, for example in the second column we have \( X_{12}:T_{21}X_1=X_2\)). Items with one index are terms of some type at a vertex, and items with two indices are terms of a type on an edge.}
\label{fig:vector_transport}
\end{figure}

As we traverse an edge, say \( e_{12} \), we get a path in \( T_2 \) which is the image of \( e_{12} \) under \( X \), denoted \( X_{12} \). As we traverse an additional edge, \( X_{12} \) is simply mapped to the next vertex by transport. The image is carried first to \( T_{32}(X_{12}) \) then to \( T_{13}\circ T_{32}(X_{12}) \).

We wish to simplify expressions such as \( T_{13}T_{32}X_{12}\cdot T_{13}X_{23}\cdot X_{31} \), which take place in a particular fiber (\( T_1 \) in this case), and which depend on arranging for the endpoint of one segment to agree with the start of another. The simplification will empower us to easily perform calculations over the whole manifold, and to prove Theorem~\ref{thm:total_index_total_curvature}.

Recall that \( \so \) acts on all the fibers of \( T \), and we have equivalences \( (\pr_1,\phi):T_i\times S^1\simto T_i\times T_i \) where \( \phi \) is the action.

\begin{mydef}
If \( (\pr_1,s) \) is the inverse of \( (\pr_1,\phi) \), the map \( s:T_i\times T_i\to S^1 \) is called \defemph{subtraction}. It maps \( (x,y) \) to the unique term \( \delta:\so \) such that \( \phi(x, \delta)=y \). We denote \( s(x,y) \) by \( y-x \).
\end{mydef}

\begin{mylemma}
\label{lem:multpath}
If \( G \) is a higher group with multiplication \( \mu:G\times G\to G \) and proof of commutativity \( \mathsf{is\underscore comm}:\pit{a,b:G}\mu(a, b)=\mu(b, a) \) then \( \mu \) induces a function \( \mu_=:(x=_G y)\times (x'=_G y')\to (\mu(x, x')=_G\mu(y,y')) \).
\end{mylemma}
\begin{myproof}
If \( p:x=_G y \) and \( p':x'=_G y' \), then we can define \( \mu_=(p, p') \) by concatenating the three paths 
\begin{align*} 
\mu(x',p)&:\mu(x', x)=_G\mu(x', y)\\
\mathsf{is\underscore comm}(x',y)&:\mu(x',y)=_G\mu(y,x')\\
\mu(y, p')&:\mu(y, x')=_G\mu(y, y').\qedhere
\end{align*}
\end{myproof}

Each fiber \( T_i \) is pointed by \( X_i \), so we can define the map \( T_i\xrightarrow[]{-X_i}\so \), and then give a name to the special term \( T_{ij}X_j - X_i \).

\begin{mydef}
We define the \defemph{rotation} of \( T_{ij} \) by \( \rho_{ij}\defeq T_{ij}X_j - X_i \). The map \( +\rho_{ij}:\so\to\so \) makes the following diagram commute%
\end{mydef}%
\begin{equation}
\begin{tikzcd}[ampersand replacement=\&]
  {T_i} \& {T_j} \\
  {S^1} \& {S^1}
  \arrow["{T_{ji}}", from=1-1, to=1-2]
  \arrow["{\alpha_i}"', from=1-1, to=2-1]
  \arrow["{\alpha_j}", from=1-2, to=2-2]
  \arrow["{\mathsf{base}\mapsto X_i}", shift left=3, curve={height=-6pt}, from=2-1, to=1-1]
  \arrow["{(-)+\rho_{ji}}"', from=2-1, to=2-2]
  \arrow["{\mathsf{base}\mapsto X_j}"', shift right=3, curve={height=6pt}, from=2-2, to=1-2]
\end{tikzcd}
\label{eq:action_of_rho}
\end{equation}

\begin{mylemma} \( \rho_{ij}+\rho_{ji}=_{\so}0 \).
\label{lem:subtraction}\end{mylemma}
\begin{myproof}
\begin{align*}
(X_i + \rho_{ij})+\rho_{ji} &=_{T_i} T_{ij}X_j + \rho_{ji}&&\text{by definition of }\rho_{ij}\\
&=_{T_i} T_{ij}(X_j + \rho_{ji})&&\text{by equivariance of transport} \\
&=_{T_i} T_{ij}(T_{ji}(X_i))&&\text{by definition of }\rho_{ji}\\
&=_{T_i} X_i&&\text{by definition of }T\qedhere
\end{align*}
\end{myproof}
We can obtain a path in \( \so \) from a dependent path by again subtracting the basepoint:
\begin{align*}
\sigma_{ij}\defeq X_{ij}-X_j&:T_{ji}X_i-X_j=_{\so}0 \\
&:\rho_{ji}=_{\so}0.
\end{align*}
Notice the reversal in indices between \( \sigma_{ij} \) and \( \rho_{ji} \), which reflects our opposite conventions for \( X_{ij} \) which we view as a path, and \( T_{ji} \) which we view as a function.

The key technical lemma is the following. Recall that dependent functions such as our vector field \( X \) send \( \refl_{v_i} \) to \( \refl_{X_i} \) (by path induction, see for example \cite{hottbook} Lemma 2.3.4). As stated, this can only be used when we have a path that is \( \refl \), for example when traversing \( e_{ij}\cdot e_{ji} \), i.e. an edge followed by its inverse. We will use subtraction together with the operation of Lemma~\ref{lem:multpath} to lift this requirement, by obtaining paths in \( \so \) that can be added without needing to be concatenated directly.

\begin{mylemma}
\label{lem:cancellation}
With the notation \( \sigma_{ij}\defeq X_{ij}-X_j \), and with addition of paths as in Lemma~\ref{lem:multpath}, we have \( \sigma_{ij}+\sigma_{ji}=_{0=_{\so}0}\refl_0 \).
\end{mylemma}
\begin{myproof}
First we need to show that the sum is a loop, then we can prove that it is \( \refl_0 \). The terms have these types:
\begin{align*}
\sigma_{ij}&:\rho_{ji}=_{\so}0\\
\sigma_{ji}&:\rho_{ij}=_{\so}0
\end{align*}
so when we add these paths with \( + \) we obtain a path \( \sigma_{ij}+\sigma_{ji}:\rho_{ji}+\rho_{ij}=_{\so} 0 \) which by Lemma~\ref{lem:subtraction} is \( 0=_{\so}0 \). We compute \( + \) using the concatenations in Lemma~\ref{lem:multpath}, which gives
\begin{align*}
\sigma_{ij}+\sigma_{ji} &: 
  (0+0) \overeq{\mathrm{Lemma~}\ref{lem:subtraction}}
  (\rho_{ji}+\rho_{ij}) \overeq{\sigma_{ij}+\rho_{ij}}
  (0+\rho_{ij}) &&\overeq{\mathsf{is\underscore comm}}
  (\rho_{ij}+0) \overeq{\sigma_{ji}}
  (0+0)\\
\sigma_{ij}+\sigma_{ji}&=_{0=_{\so}0} (\sigma_{ij}+\rho_{ij})\cdot \sigma_{ji}&& \\
&=_{0=_{\so}0} ((X_{ij}-X_j)+\rho_{ij})\cdot (X_{ji} - X_i)&&\text{definition of }\sigma s\\
&=_{0=_{\so}0} (T_{ij}(X_{ij})-X_i)\cdot (X_{ji} - X_i)&&\text{action of }\rho_{ij}\text{ (see (\ref{eq:action_of_rho}))} \\
&=_{0=_{\so}0} \refl_{X_i}-X_i&&\text{path induction for }X\\
&=_{0=_{\so}0} \refl_0&&\qedhere
\end{align*}
\end{myproof}
\begin{mynote}
The classical argument of Hopf \cite{hopf}, which is presented in more detail in the more readily available \cite{needham}, makes an implicit assumption that we can concatenate two terms of different type such as \( X_{ij} \) and \( X_{kl} \). The authors name such terms ``change in angle across an edge'' and stipulate that this is a function ``defined on edges.'' The extra work we are doing in this section amounts to a partial formalization of this idea.
\end{mynote}

\subsection{An example vector field on the sphere}
\label{sec:octahedron_vector_field}
\begin{figure}[h]
\centering
\begin{tikzpicture}%
  [x={(-0.860769cm, -0.121512cm)},
  y={(0.508996cm, -0.205391cm)},
  z={(-0.000053cm, 0.971107cm)},
  scale=1,
  back/.style={loosely dotted, thin},
  edge/.style={black, thick},
  arrow/.style={black, very thick, solid, -{Stealth[scale=0.8]}},
  facet/.style={fill=blue!95!black,fill opacity=0.0},
  vertex/.style={inner sep=1pt,circle,draw=green!25!black,fill=black,thick}]
\begin{scope}[nodes=vertex]
\node[label=above right:\( b \)] at (-1, 1, 0) (b)     {};
\node[label=below:\( y \)] at (0, 0, -1.4) (y)    {};
\node[label=above:\( w \)] at (0, 0, 1.4)  (w)   {};
\node[label=above left:\( g \)] at (1, -1, 0) (g)    {};
\node[label=above left:\( r \)] at (1, 1, 0)  (r)   {};
\node[label=above right:\( o \)] at (-1, -1, 0) (o)    {};
\end{scope}
\draw[edge,back,arrow] (o) -- (b);
\draw[edge,back,arrow] (y) -- (o);
\draw[edge,back] (o) -- (w);
\draw[edge,back] (o) -- (g);
\node[vertex] at (o)     {};
\fill[facet] (1, 1, 0) -- (0, 0, -1.4) -- (1, -1, 0) -- cycle {};
\fill[facet] (1, 1, 0) -- (0, 0, 1.4) -- (1, -1, 0) -- cycle {};
\fill[facet] (1, 1, 0) -- (-1, 1, 0) -- (0, 0, 1.4) -- cycle {};
\fill[facet] (1, 1, 0) -- (-1, 1, 0) -- (0, 0, -1.4) -- cycle {};
\draw[edge,arrow] (b) -- (y);
\draw[edge] (b) -- (w);
\draw[edge] (b) -- (r);
\draw[edge] (y) -- (g);
\draw[edge] (y) -- (r);
\draw[edge,arrow] (g) -- (w);
\draw[edge,arrow] (w) -- (r);
\draw[edge,arrow] (r) -- (g);
\end{tikzpicture}
\caption{A vector field on \( \oo \), which extends to a rotation of the octahedron in space.}
\label{fig:octahedron_vector_field}
\end{figure}
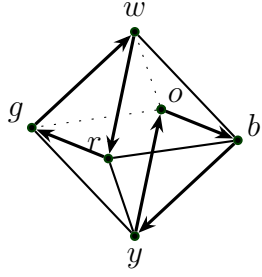

Figure~\ref{fig:octahedron_vector_field} shows an example of a vector field \( \xso \) (\( s \) for ``spin'') on the octahedron \( \oo \). The picture can really only convey the value of \( \xso \) on vertices, which it does by displaying it as an arrow on the surface itself, for example as an arrow from \( w \) to \( r \), instead of trying to draw the fiber at \( w \) which is where \( r:\link(w) \) actually lives. But hopefully there is a net gain in clarity.

We will define \( \xso \) on edges, and then compute the swirling around each face. These calculations will then be cited later when we discuss the total swirling on a realization. Here is all the data for \( \xso \), which includes values when traversing an edge in either direction:
\[\begin{array}{|c|c|}
\hline
\multicolumn{2}{|c|}{\xso\text{ on vertices}} \\
\hline
v & \xso(v) \\
\hline
w & r \\
r & g \\
g & w \\
\hline
o & b \\
b & y \\
y & o \\
\hline
\end{array}\,
\begin{array}{|c|c|}
\hline
\multicolumn{2}{|c|}{\xso\text{ on north edges}} \\
\hline
e & \xso(e) \\
\hline
wr & yg \\
wg & rw \\
wo & wb \\
wb & ry \\
\hline
rw & gr \\
gw & br \\
ow & br \\
bw & br \\
\hline
\end{array}\,
\begin{array}{|c|c|}
\hline
\multicolumn{2}{|c|}{\xso\text{ on south edges}} \\
\hline
e & \xso(e) \\
\hline
yr & yg \\
yg & ow \\
yo & wb \\
yb & oy \\
\hline
ry & go \\
gy & go \\
oy & bo \\
by & go \\
\hline
\end{array}\,
\begin{array}{|c|c|}
\hline
\multicolumn{2}{|c|}{\xso\text{ on equator edges}} \\
\hline
e & \xso(e) \\
\hline
br & yg \\
rg & ow \\
go & wb \\
ob & ry \\
\hline
rb & ry \\
gr & wg \\
og & ow \\
bo & yb \\
\hline
\end{array}\]
How do we read the tables for the edges? For example, \( \xso(wr) \) should be a term of type \( \tr(wr)\xso(w)=_{Tr}\xso(r) \). Looking up \( \xso \) on vertices means this type is \( \tr(wr)(r)=_{Tr}g \), and looking up \( \tr \) in Definition~\ref{def:transport_oo} means the type is \( y=_{Tr}g \). We can choose \( \xso(wr) \) to be any path from \( y \) to \( g \) and we choose the short path \( yg \).

Next we want to compute \( \xso \) around a face. We will do this for the four north faces only. In this calculation we will denote paths by listing the vertices, so \( rg\cdot gw \) is denoted \( rgw \) and \( wrgw \) is a loop at \( w \).
\begin{align*}
\xso(wrgw) &= \tr(rgw)\xso(wr)\cdot\tr(gw)\xso(rg)\cdot\xso(gw)\\
&= \tr(rgw)(yg)\cdot\tr(gw)(ow)\cdot(br)\\
&= (yg\mapsto yo\mapsto \boxed{go})\cdot(ow\mapsto \boxed{ob})\cdot \boxed{br}\\
&= \boxed{gobr}
\end{align*}
To help follow along, in the third line we are showing the results of intermediate transports, and put a box around the result. The other three faces give
\begin{align*}
\xso(wgow) &= \tr(gow)\xso(wg)\cdot\tr(ow)\xso(go)\cdot\xso(ow)\\
&= \tr(gow)(rw)\cdot\tr(ow)(wb)\cdot(br)\\
&= (rw\mapsto gw\mapsto \boxed{gr})\cdot(wb\mapsto \boxed{rb})\cdot \boxed{br}\\
&= \boxed{gr}
\end{align*}

\begin{align*}
\xso(wobw) &= \tr(obw)\xso(wo)\cdot\tr(bw)\xso(ob)\cdot\xso(bw)\\
&= \tr(obw)(wb)\cdot\tr(bw)(ry)\cdot(br)\\
&= (wb\mapsto wr\mapsto \boxed{gr})\cdot(ry\mapsto \boxed{rb})\cdot \boxed{br}\\
&= \boxed{gr}
\end{align*}

\begin{align*}
\xso(wbrw) &= \tr(brw)\xso(wb)\cdot\tr(rw)\xso(br)\cdot\xso(rw)\\
&= \tr(brw)(ry)\cdot\tr(rw)(yg)\cdot(gr)\\
&= (ry\mapsto gy\mapsto \boxed{gr})\cdot(yg\mapsto \boxed{rg})\cdot \boxed{gr}\\
&= \boxed{gr}
\end{align*}

So we have these four edges in \( Tw \), which is pointed at \( \xso(w)=r \): \( gobr, gr, gr, gr \). Lemma~\ref{lem:multpath} tells us we can add these by transporting them so as to be concatenate-able. This gives \( gobr\cdot rb\cdot bo\cdot og=\refl_g \) in \( Tw \). So the total swirling in the northern hemisphere is trivial. We leave the computation in the southern hemisphere as an exercise, but it is very symmetrical with the northern hemisphere, and gives \( \refl_o \) in \( Ty \).

\clearpage
\section{The total construction}
\label{sec:totals}
We will place holonomy, flatness, and vector fields on the same footing, and combine them. We will prove a relation between the total curvature of a tangent bundle, and total index of a vector field. This is the key step in the simultaneous proof of the Gauss-Bonnet theorem and the Poincaré-Hopf theorem, lacking only the relationship with the Euler characteristic.

\subsection{Index of the vector field on a face}
The index of a vector field is derived from other data that we have. Consider a 2-dimensional simplicial complex \( M \) and its 2-dimensional realization \( \mm=\mm_0\to\mm_1\to\mm_2 \), and a principal circle bundle \( T:\mm\to\Kzt \) extending \( \link \). Consider a face \( F:M_2 \), with \( m:\mm \) a vertex of \( F \) and \( \partial F:m=m \) the boundary loop of the face. We have accumulated the following constructions:
\[\begin{aligned}
\tr_F&\defeq \tr(\partial F)&&:Tm=Tm&&\text{holonomy}\\
\flat_F&\defeq \flat(\partial F)&&:\id=_{Tm=Tm}\tr(\partial F)&&\text{flatness}\\
X_F&\defeq X(\partial F)&&:\tr(\partial F)(X(m))=_{Tm}X(m)&&\text{swirling.}\\
\end{aligned}\]
The flatness is a path in automorphisms of \( Tm \) whereas the swirling is a path in \( Tm \), but we can view the latter as a path in automorphisms as well:
\begin{myprop}
\label{prop:eveq}
Given a polygon \( \ccc(n):\Kzt \) and a point \( b:\ccc(n) \) the evaluation map \( \ev(b):(\ccc(n)=_{\Kzt}\ccc(n))\to \ccc(n) \) is an equivalence.
\end{myprop}
\begin{myproof}
See \cite{buchholtz2023central}.
\end{myproof}

Now \( \sigma^X_F\defeq\ev(X(m))^{-1}(X_F) \) is the automorphism of \( Tm \) that corresponds to the swirling. We now have:
\[\begin{aligned}
\tr_F&:Tm=Tm&&\text{holonomy on }F\\
\flat_F&:\id=_{Tm=Tm}\tr_F&&\text{flatness on }F\\
\sigma^X_F&:\tr_F=_{Tm=Tm}\id&&\text{swirling on }F.\\
\end{aligned}\]
These last two can be concatenated to make a loop. We can then obtain an integer by noting that the automorphisms are a circle \( \loopy_{Tm}(\Kzt)\simeq S^1 \), and recalling the well known fact that loops in the circle are integers: \( \loopy(S^1,\base)\simeq\zz \) (e.g. \cite{hottbook} Corollary 8.1.10).
\begin{mydef}
\label{def:index}
The \defemph{index of the vector field \( X \) on the face \( F \)} is the integer \( I^X_F\defeq\flat_F\cdot \sigma^X_F:(\id=_{Tm=Tm}\id) \).
\end{mydef}
We now have the final list of ingredients at a single face:
\begin{equation}
\label{eq:face_elements}
\begin{aligned}
\tr_F&:Tm=Tm&&\text{holonomy on }F\\
\flat_F&:\id=_{Tm=Tm}\tr_F&&\text{flatness on }F\\
\sigma^X_F&:\tr_F=_{Tm=Tm}\id&&\text{swirling on }F\\
I^X_F&:\zz&&\text{index on }F.
\end{aligned}
\end{equation}

\subsection{Total flatness, total index}
Now we wish to compute a sum over all faces of the data. To do this we need all of our definitions plus two new ones that help us organize the sum and make a key assumption about cancellation of edges.

\begin{mydef}
A 2-dimensional simplicial complex \( M^+=[M_0, M_1, M_2] \) with orientation \( \mathscr{O} \) is \defemph{without boundary} if every edge is a subset of exactly two faces, and if the orderings induced by the two inclusions differ by a swap of vertices.
\end{mydef}

Some of the assumptions in these definitions are unnecessary because they are entailed in other assumptions, but we will not attempt to eliminate them.

\begin{mydef}
An \defemph{oriented set of pointed faces} of a realization \( \mm \) of a 2-dimensional simplicial complex \( M=[M_0,M_1,M_2] \) with orientation \( \mathscr{O} \) is a choice \( v_F:F \) of a vertex in each face \( F:M_2 \), and the boundary path \( \ell_F:v_F=v_F \) that concatenates the three edges of the face in the order determined by the orientation.
\end{mydef}

For the rest of the note we will equip our simplicial complex \( M \) with an orientation \( \mathscr{O} \), and our realization \( \mm \) with an oriented set of pointed faces.

We have seen in Lemma~\ref{lem:multpath} how to perform the required sum, namely by using subtraction to map all \( T_{v_F} \) to \( S^1 \) by subtracting \( X(v_F) \).
\begin{mydef}
The \defemph{net holonomy}, \defemph{total flatness}, \defemph{total swirling}, and \defemph{total index} on an oriented set of pointed faces is
\begin{equation}
\label{eq:tot_elements}
\begin{aligned}
\tr_\tot      & \defeq \sum_F\tr_F &&:\so\\
\flat_\tot    & \defeq \sum_F\flat_F &&:\base=_{\so}\tr_\tot\\
\sigma^X_\tot & \defeq \sum_F\sigma^X_F &&:\tr_\tot=_{\so}\base\\
I^X_\tot      & \defeq \sum_F I^X_F &&:\zz.
\end{aligned}
\end{equation}
\end{mydef}

The assumption about edges appearing twice, once in each direction, together with Lemma \ref{lem:subtraction} proves the following
\begin{myprop}The net holonomy is \( 0 \) in \( \so \).\label{prop:total_rotation}
\end{myprop}
\begin{myproof}
With a vector field \( X \) chosen, we can compute \[ \sum_F\tr_F=\left[\sum_{(\mathrm{edges\ }e_{ij})} \rho_{ij}+\rho_{ji}\right]=_{\so}0 \] and note that the result does not depend on \( X \).
\end{myproof}

\begin{mycor}
Total flatness is a loop: \( \flat_\tot:\base=_{\so}\base \).
\end{mycor}\begin{myproof}Follows from Corollary~\ref{prop:total_rotation}.\end{myproof}

And making use of Lemma \ref{lem:cancellation} we obtain
\begin{myprop}The total swirling is \( \refl_0 \) in \( \so \).\label{prop:total_swirling}
\end{myprop}
\begin{myproof}We can compute
\[\sum_F\sigma^X_F=\left[\sum_{(\mathrm{edges\ }e_{ij})}\sigma_{ij}+\sigma_{ji}\right]=_{\so}\refl_0.\qedhere\]
\end{myproof}

And finally, if we refer to the passage from a loop \( \ell:\base=_{\so} \base \) to the corresponding integer in \( \loopy\so \) as ``winding number'':

\begin{mythm}\label{thm:total_index_total_curvature}The winding number of total flatness equals the total index.
\end{mythm}
\begin{myproof}
Follows directly from \( I^X_\tot\defeq\flat_\tot\cdot \sigma^X_\tot \) and Proposition~\ref{prop:total_swirling}.
\end{myproof}

\begin{mynote}
The Gauss-Bonnet theorem is often stated in the form \[\frac{1}{2\pi}\int_M K\,dA=\chi(M)\] and although we won't attempt an exhaustive translation between this and our approach, we want to point out that the factor of \( \frac{1}{2\pi} \) is doing exactly what winding number is doing for us. It seems likely that attending fully to constant coefficients in other important formulas can yield conceptual fruit.
\end{mynote}

\subsection{Calculation on the sphere}
The flatness structures on the eight faces of \( \oo \) in Definition~\ref{def:octahedron_curvature} are all identical and equal to \( H_R:\id=_{\ccc(4)=\ccc(4)}R \) on a fixed 4-gon. Composing \( R \) 8 times gives the identity, as required by Proposition~\ref{prop:total_rotation}. Adding the flatness structures gives \( 8H_R:\id=_{\ccc(4)=\ccc(4)}\id \) whose winding number is 2.

Next we want to compute the total index of the vector field \( \xso \) from Section~\ref{sec:octahedron_vector_field}. We saw that on the four faces in the northern hemisphere the swirling gave us these four paths in \( Tw \): \( gobr, gr, gr, gr \). The flatness structures on each of these faces is \( H_R:\id=_{Tw=Tw}R \) where \( R \) rotates by one vertex clockwise. Using Definition~\ref{def:index} for the index we see that we first need to convert \( gobr \) and \( gr \) into paths of rotations, which give
\begin{align*}
\sigma^{\xso}_{wrgw}=H_R^3   &:R^{-3}=id \\
\sigma^{\xso}_{wgow}=H_R^{-1}&:R=id \\
\sigma^{\xso}_{wobw}=H_R^{-1}&:R=id \\
\sigma^{\xso}_{wbrw}=H_R^{-1}&:R=id \\
\end{align*}
Now forming the concatenation of \( H_R \) with each of these gives \( H_R^4, \refl_\id, \refl_\id, \refl_\id \) respectively. The winding numbers are \( 1, 0, 0, 0 \). So the index around the face that has all the swirling in Figure~\ref{fig:octahedron_vector_field} is 1 and around the other three is 0. We haven't written out the computation for the southern hemisphere but it will provide another 1 and three 0s, giving a total index of 2.

\clearpage
\appendix
\section{Appendix: the program ahead}

Follow-up work could remove the need to assume that maps from the realization of oriented simplicial complexes factor through \( \Kzt \), instead proving that this is the case. This would unify two assumptions of orientability into one. The assumption that edges appear twice, once with each orientation, could also be deduced from presumably standard material.

An internal definition of Euler characteristic is needed to provide the relationship between total index and the classic invariant.

More relationships with discrete differential geometry\cite{crane_ddg}\cite{crane_connections} could be explored.

The results could be formalized. A good starting point is Lemma~\ref{lem:octahedron_sphere} that the octahedron is equivalent to the 2-sphere. A good next step would be to prove that a combinatorial \( n \)-sphere has realization equivalent to the HoTT \( n \)-sphere. This would internalize the function \( \link \).

Our hope is that this note can form the starting point for a complete formulation of the homotopical results of gauge theory and Chern-Weil theory. Classically we know the space of connections is contractible, and the group of based gauge transformations (see below) acts freely on it. Those statements should be easy to make in our framework. In \cite{atiyah1983yang} connections are presented as a splitting of a short exact sequence of bundles. The community has grappled with this picture (see the nLab at \cite{urs_atiyah}), and it deserves clarification in HoTT.

The paper of Freed and Hopkins\cite{freed2013chernweil} seeks a classifying space for the space of connections on a principal bundle. Their paper served as a primary motivation for this work, and we hope we have taken a step in that direction. Its focus on homotopy theory, model structures, and groupoids should make it accessible to many HoTT theoreticians. We conjecture that the theory of Koszul complexes that arises is the infinitesimal shadow of the higher groupoid structure of the relevant HoTT classifying space.

\subsection{Gauge transformations}
\begin{mydef}
\label{sec:automorphisms}
Suppose we have \( T:M\to\EMzo \) and \( P\defeq\sit{x:M}Tx \). Then we can form the type family \( \Aut T:M\to\uni \) given by \( \Aut T(x)\defeq(Tx=Tx) \). The total space \( \Aut P\defeq\sit{x:M}(Tx=Tx) \), which is a bundle of groups, is called the \defemph{automorphism bundle} or the \defemph{gauge bundle} and sections \( \pit{x:M}(Tx=Tx) \), which are homotopies \( T\sim T \), are called \defemph{automorphisms of \( P \)} or \defemph{gauge transformations}. If \( m:M \) is a basepoint and if \( p:Tm \) is a basepoint in the fiber at \( m \), then the \defemph{pointed automorphisms of \( P \)} or \defemph{based gauge transformations}.
\end{mydef}

\clearpage

\bibliography{connections}

\begin{thebibliography}{10}
\providecommand{\url}[1]{#1}
\csname url@samestyle\endcsname
\providecommand{\newblock}{\relax}
\providecommand{\bibinfo}[2]{#2}
\providecommand{\BIBentrySTDinterwordspacing}{\spaceskip=0pt\relax}
\providecommand{\BIBentryALTinterwordstretchfactor}{4}
\providecommand{\BIBentryALTinterwordspacing}{\spaceskip=\fontdimen2\font plus
\BIBentryALTinterwordstretchfactor\fontdimen3\font minus
  \fontdimen4\font\relax}
\providecommand{\BIBforeignlanguage}[2]{{%
\expandafter\ifx\csname l@#1\endcsname\relax
\typeout{** WARNING: IEEEtran.bst: No hyphenation pattern has been}%
\typeout{** loaded for the language `#1'. Using the pattern for}%
\typeout{** the default language instead.}%
\else
\language=\csname l@#1\endcsname
\fi
#2}}
\providecommand{\BIBdecl}{\relax}
\BIBdecl

\bibitem{talbott}
\BIBentryALTinterwordspacing
S.~Talbott, \emph{The Future Does Not Compute: Transcending the Machines in Our
  Midst}.\hskip 1em plus 0.5em minus 0.4em\relax O'Reilly \& Associates, 1995.
  [Online]. Available: \url{https://books.google.com/books?id=KcXaAAAAMAAJ}
\BIBentrySTDinterwordspacing

\bibitem{hopf}
\BIBentryALTinterwordspacing
H.~Hopf, ``Differential geometry in the large: Seminar lectures, new york
  university, 1946 and stanford university, 1956,'' 1983. [Online]. Available:
  \url{https://api.semanticscholar.org/CorpusID:117042538}
\BIBentrySTDinterwordspacing

\bibitem{needham}
\BIBentryALTinterwordspacing
T.~Needham, \emph{Visual Differential Geometry and Forms: A Mathematical Drama
  in Five Acts}.\hskip 1em plus 0.5em minus 0.4em\relax Princeton University
  Press, 2021. [Online]. Available:
  \url{https://books.google.com/books?id=Mc0QEAAAQBAJ}
\BIBentrySTDinterwordspacing

\bibitem{ljungstrom_cellular}
A.~Ljungström and L.~Pujet, ``Cellular methods in homotopy type theory,''
  \url{https://pujet.fr/pdf/cellular.pdf}, 2025.

\bibitem{buchholtz2023central}
U.~Buchholtz, J.~D. Christensen, J.~G.~T. Flaten, and E.~Rijke, ``Central
  h-spaces and banded types,'' 2023.

\bibitem{sco}
\BIBentryALTinterwordspacing
L.~Scoccola, ``Nilpotent types and fracture squares in homotopy type theory,''
  \emph{Mathematical Structures in Computer Science}, vol.~30, no.~5, p.
  511–544, May 2020. [Online]. Available:
  \url{http://dx.doi.org/10.1017/S0960129520000146}
\BIBentrySTDinterwordspacing

\bibitem{hottbook}
{Univalent Foundations Program}, \emph{Homotopy Type Theory: Univalent
  Foundations of Mathematics}.\hskip 1em plus 0.5em minus 0.4em\relax Institute
  for Advanced Study: \url{https://homotopytypetheory.org/book}, 2013.

\bibitem{dcct}
\BIBentryALTinterwordspacing
U.~Schreiber, ``\BIBforeignlanguage{english}{Differential cohomology in a
  cohesive infinity-topos},'' Oct. 2013. [Online]. Available:
  \url{https://arxiv.org/abs/1310.7930v1}
\BIBentrySTDinterwordspacing

\bibitem{myersgood}
D.~J. Myers, ``Good fibrations through the modal prism,'' 2022.

\bibitem{egbert}
E.~Rijke, \emph{Introduction to Homotopy Type Theory}, ser. Cambridge Studies
  in Advanced Mathematics.\hskip 1em plus 0.5em minus 0.4em\relax Cambridge
  University Press, 2025.

\bibitem{whitehead_triangulation}
J.~H.~C. Whitehead, ``On {$C^1$}-complexes,'' \emph{Annals of Mathematics}, pp.
  809--824, 1940.

\bibitem{kirby_siebenmann}
R.~C. Kirby and L.~C. Siebenmann, \emph{Foundational essays on topological
  manifolds, smoothings, and triangulations}, ser. Annals of Mathematics
  Studies.\hskip 1em plus 0.5em minus 0.4em\relax Princeton University Press,
  1977, no.~88, with notes by John Milnor and Michael Atiyah. MR:0645390.
  Zbl:0361.57004.

\bibitem{kobayashinomizu}
\BIBentryALTinterwordspacing
S.~Kobayashi and K.~Nomizu, \emph{Foundations of Differential Geometry}, ser.
  Foundations of Differential Geometry.\hskip 1em plus 0.5em minus 0.4em\relax
  Interscience Publishers, 1963. [Online]. Available:
  \url{https://books.google.com/books?id=wn4pAQAAMAAJ}
\BIBentrySTDinterwordspacing

\bibitem{crane_ddg}
\BIBentryALTinterwordspacing
K.~Crane, F.~de~Goes, M.~Desbrun, and P.~Schröder, ``Digital geometry
  processing with discrete exterior calculus,'' in \emph{ACM SIGGRAPH 2013
  courses}, ser. SIGGRAPH '13.\hskip 1em plus 0.5em minus 0.4em\relax New York,
  NY, USA: ACM, 2013. [Online]. Available:
  \url{{https://www.cs.cmu.edu/~kmcrane/Projects/DDG/}}
\BIBentrySTDinterwordspacing

\bibitem{crane_connections}
\BIBentryALTinterwordspacing
K.~Crane, ``Discrete connections for geometry processing,'' Master's thesis,
  California Institute of Technology, 2010. [Online]. Available:
  \url{http://resolver.caltech.edu/CaltechTHESIS:05282010-102307125}
\BIBentrySTDinterwordspacing

\bibitem{atiyah1983yang}
M.~F. Atiyah and R.~Bott, ``The yang-mills equations over riemann surfaces,''
  \emph{Philosophical Transactions of the Royal Society of London. Series A,
  Mathematical and Physical Sciences}, vol. 308, no. 1505, pp. 523--615, 1983.

\bibitem{urs_atiyah}
{nLab authors}, ``{{A}}tiyah {{L}}ie groupoid,''
  \url{https://ncatlab.org/nlab/show/Atiyah+Lie+groupoid}, Oct. 2024,
  \href{https://ncatlab.org/nlab/revision/Atiyah+Lie+groupoid/25}{Revision 25}.

\bibitem{freed2013chernweil}
\BIBentryALTinterwordspacing
D.~S. Freed and M.~J. Hopkins, ``Chern-weil forms and abstract homotopy
  theory,'' 2013. [Online]. Available: \url{https://arxiv.org/abs/1301.5959}
\BIBentrySTDinterwordspacing

\end{thebibliography}
\end{document}